\documentclass[1p]{elsarticle}

\journal{Applied Mathematics and Computation}

\usepackage{amsmath}
\usepackage{commath}
\usepackage{amssymb}
\usepackage{amsthm}
\usepackage{hyperref}
\usepackage{textcomp}

\usepackage{float}
\usepackage{booktabs}

\newtheorem{theorem}{Theorem}

\newtheorem{criterion}{Criterion}

\newtheorem{remark}{Remark}
\newtheorem{definition}{Definition}
\newcommand{\R}{\mathbb{R}}

\usepackage[dvipsnames]{xcolor}

\usepackage{etoolbox}
\patchcmd{\pprintMaketitle}
  {\fi\hrule}
  {\fi\ifvoid\extrainfobox\else\unvbox\extrainfobox\par\vskip10pt\fi\hrule}
  {}{}

\newenvironment{extrainfo}
  {\global\setbox\extrainfobox=\vbox\bgroup\parindent=0pt }
  {\egroup}
\newsavebox\extrainfobox

\begin{document}

\begin{frontmatter}
\title{On the choice of initial guesses \\for the Newton-Raphson algorithm}
\author{Francesco Casella\corref{cor1}}
\ead{francesco.casella@polimi.it}
\cortext[cor1]{Corresponding author}
\address{Dipartimento di Elettronica, Informazione e Bioingegneria\\
        Politecnico di Milano, Italy}
\author{Bernhard Bachmann}
\ead{bernhard.bachmann@fh-bielefeld.de}
\address{Faculty of Engineering and Mathematics\\
         University of Applied Sciences Bielefeld, Germany}
\fntext[]{Paper accepted for publication, see \url{https://doi.org/10.1016/j.amc.2021.125991}}
\fntext[]{(C) 2021. This manuscript version is made available under the CC-BY-NC-ND 4.0 license \url{http://creativecommons.org/licenses/by-nc-nd/4.0/}}
\begin{abstract}
The initialization of equation-based differential-algebraic system models, and more in general the solution of many engineering and scientific problems, require the solution of systems of nonlinear equations. Newton-Raphson's method is widely used for this purpose; it is very efficient in the computation of the solution if the initial guess is close enough to it, but it can fail otherwise. In this paper, several criteria are introduced to analyze the influence of the initial guess on the evolution of 
Newton-Raphson's algorithm and to identify which initial guesses need to be improved in case of convergence failure. In particular, indicators based on first and second derivatives of the residual function are introduced, whose values allow to assess how much the initial guess of each variable can be responsible for the convergence failure. The use of such criteria, which are based on rigorously proven results, is successfully demonstrated in three exemplary test cases.
\end{abstract}

\begin{keyword}
Newton-Raphson's algorithm, Convergence, Nonlinear equations, Equation-based modelling.
\end{keyword}

\begin{extrainfo}
Declarations of interest: none.
\end{extrainfo}

\end{frontmatter}

\section{Introduction} \label{sec:Introduction and Motivation}
\subsection{Goal of the paper}
\label{sec:Goal}
Newton-Raphson's (NR) algorithm and its variants have been used for over 250
years to solve implicit nonlinear equations. The algorithm is iterative and the
convergence to the desired solution crucially depends on the choice of the
initial guess for the unknowns of the problem. Once the result of the iterations
is close enough to the solution, under mild regularity conditions and under the
assumption of non-singular Jacobian, the algorithm converges to the solution in
a superlinear fashion.

In general, it may not be easy or practical to obtain an initial guess close
enough to the solution to ensure that the asymptotic convergence result is
obtained after a small number of iterations. In fact, if the initial guess is
sufficiently far from the sought-after solution, NR's algorithm may not converge
at all to it.

The asymptotic convergence properties of NR's algorithm when the initial guess is close enough to the solution are well-known, see e.g. \cite{Ortega1970}. 

Many theorems have been proven in the past, originated from the classical result
by Kantorovich \cite{Kantorovich1964}, which provide sufficient conditions for
the convergence of NR's method to a solution, see e.g. \cite{Galantai2000,Argyros2012} and references therein. Although these results are very powerful from a mathematical point of view, in most practical cases it is hardly possible to verify the required convergence conditions. 

Problem-specific methods are proposed in the literature to improve the chances and speed of convergence of NR's method. For example, \cite{VolmerEtAl2020} proposes a method to improve the convergence of the solution of visco-plastic models; \cite{CherifEtAl2019} proposes a strategy to obtain good initial guesses for magnetostatic problems; \cite{SunEtAl2017} proposes a convergence analysis for the electrical power flow problem, which allows to obtain initial guesses close to the solution. However, each of these methods is limited to a very specific class of mathematical models and cannot be extended to handle generic systems of equations, possibly describing multi-domain systems.

Up to the authors' knowledge, the following two questions have not received an answer so far in the published literature:
\begin{itemize}
\item For which unknowns of a nonlinear system is it actually necessary to provide a good initial guess to NR's method?
\item In case of convergence failure of NR's method starting from a certain initial guess, how should one improve it to eventually achieve convergence to the desired solution?
\end{itemize}

These questions arise in a large number of practical cases, whenever NR's method is used to solve a nonlinear problem of any kind. The aim of this paper is to answer these two questions based as much as possible on rigorous results and, where necessary, on heuristic criteria based on those rigorous results.

\subsection{Background}
The introduction of Equation-based, Object-Oriented modelling Languages and
Tools (EOOLTs), such as Modelica
\cite{MattssonElmqvistOtterCEP1998,Fritzson2014} or gPROMS
\cite{BartonPandelidesAIChE1994}, started in the mid 90's of the last century,
has made the need of good answers to the above-mentioned questions
a compelling requirement.

These modelling languages allow to build complex system models, described by
differential-algebraic equations, which can potentially span multiple physical
domains, such as mechanical, electrical, thermal, thermal-hydraulic, chemical,
etc. The system models are obtained by assembling equation-based component
models in a modular way, possibly taking them from libraries of well-tested and
validated reusable component models developed by third-parties.

The simulation of such models requires finding a consistent initial solution for
their DAEs \cite{PantelidesJSciComp1988}. Such solution can be
obtained by adding a set of initial equations to
the DAEs, often resulting in very large sets of nonlinear equations. The tools
which handle these models are agnostic, in the sense that they are not limited
to a specific physical domain or to a pre-specified set of models, for which
some heuristic criterion to provide initial guesses for the initialization
problem can be found and embedded in the software. To the contrary, the users of
such tools have complete freedom to combine equation-based models from multiple
reusable libraries, together with other models and initial conditions that they
write themselves in the form of equations with arbitrary structure.

Thanks to the high-level, modular approach to modelling of EOOLTs, building
complex system models using these tools is relatively straightforward;
unfortunately, the solution of the corresponding initialization problems often
turns out to be a critical task. In most cases, the model equations are
nonlinear, the nonlinear solvers often fail, and the end user is
left to struggle with very low-level error messages and log files to analyze, in
order to understand how to eventually succeed in the solution process. In many cases, the initialization problem is eventually solved, and no trace is left of this effort in the publications that describe the results. However, in many other
cases, this issue turns out to be the main limiting factor in the adoption of EOOLTs in a certain application domain.

There are some strategies to avoid this kind of problems. First and foremost,
whenever nonlinear equations are involved, one should have at least some idea
about the solution, and set initial guess values for the unknowns of the nonlinear
equations accordingly.

When dealing with large interconnected modular systems, one can first try to solve the
initialization problem for each of its components or sub-systems separately, applying
suitable boundary conditions, and then collect the found solutions to get the initial
guesses for the solution of the system-wide initialization process. However,
this process is quite convoluted, and it requires the \emph{a-priori} knowledge of
consistent system boundaries for the involved sub-systems.

Another well-known strategy to make the selection of initial guesses less
cumbersome (among other things) is tearing
\cite{CarpanzanoMaffezzoniMCS1998,CarpanzanoMCMDS2000,CellierKofman2006},
whereby a certain subset of tearing variables is chosen, so that all other
variables can be explicitly computed from them following a certain sequence of
assignments. In this case, it is only necessary to provide initial guesses for
the tearing variables, since guesses of all the other unknowns are automatically
computed by the sequence of assignments during the first iteration of NR's
algorithm. However, the choice of those initial guesses is still critical and
can lead to convergence failure.

In some cases, a homotopy-based approach can help finding the initial solution
of the system by first solving a simplified problem, and then transforming the
simplified problem into the actual one by means of a homotopy transformation,
using a continuation solver, see e.g.
\cite{SielemannEtAlModelica2011,SielemannCasellaOtterSIMPAT2013} and references
therein. However, also in this case, unless the simplified model is linear,
solving the initial simplified problem also requires the use of an iterative,
NR-type solver, which is prone to failure if not initialized correctly.

Summing up, despite the best efforts of the model developer, it is often the case
that some initial guess values are not close enough to the solution, or are computed
incorrectly, or possibly not set at all due to some oversight, causing NR's solver
to fail, or possibly to converge to an unwanted solution. The modeller is then left
with the task of finding out which initial guesses are not good enough and fix
them until convergence is achieved. This task is in general ugly, time-consuming,
and requiring an unknown and possibly very large amount of time, particularly
in the case of large models.

\subsection{Contents of the paper}
Given this scenario, there is a definite need of general criteria to pinpoint those
initial guesses that are causing the convergence failure of NR solvers, applied to generic
systems of nonlinear equations coming from physical system models of
arbitrary nature. To the authors' best knowledge, this general
problem is not addressed as such in the published literature. The goal of this
paper is thus to provide such criteria, based as much as possible on rigorous
results and, where necessary, on some heuristic assumptions.

The paper is structured as follows. In Section \ref{sec:Method}, several new
theorems are stated, which provide the rigorous groundwork for the remainder of
the paper. In Section \ref{sec:Discussion}, the relevance of these theorems with
respect to the two questions stated in Section \ref{sec:Goal} is discussed, leading
to the formulation of several heuristic criteria aimed at the effective
choice of initial guesses for NR's algorithm. In Section \ref{sec:Examples}, those
criteria are successfully demonstrated on three exemplary physical
modelling problems. Finally, Section \ref{sec:Conclusion} concludes the paper.

\section{Method}
\label{sec:Method}
Consider the equation
\begin{equation}
\label{eq:ImplicitEq}
f(x) = 0,
\end{equation}
where $x \in {\R}^m$ and $f:\R^m \rightarrow \R^m$ is a vector function which is
continuously differentiable in an open neighbourhood $\mathfrak{D}$ of the solution
$\bar{x}$, $f(\bar{x}) = 0$. Denote the Jacobian matrix of function $f(x)$ with respect to $x$ as $f_x(x)$. Assume the vector of the unknowns $x$ is suitably ordered,
so that it can be split into two sub-vectors $w\in\R^q$ and $z\in\R^{m-q}$
\begin{equation}
x = \begin{bmatrix} w \\ z \end{bmatrix},
\end{equation}
$w$ being the smallest possible sub-set of $x$ such that
\begin{equation}
\label{eq:JacDependency}
f_x(x) = J(w),
\end{equation}
i.e., the Jacobian matrix of $f(x)$ depends only on $w$ and not on $z$, and therefore the function $f(x)$ depends only linearly on $z$. Assume the equations in
\eqref{eq:ImplicitEq} are ordered so that $f(x)$ can be split into two vector functions
$n(x)$ and $l(x)$, $n:\R^m \rightarrow \R^p$, $l:\R^m \rightarrow \R^{m-p}$
\begin{equation}
f(x) = \begin{bmatrix} n(x) \\ l(x) \end{bmatrix},
\end{equation}
where $n(x)$ contains the non-linear equation residuals and $l(x)$ contains the linear
equation residuals.

The solution $\bar{x}$ can be computed iteratively by NR's algorithm, which requires to
solve the following linear equation at each iteration $j$
\begin{equation}
\label{eq:Newton}
f_x(x_{j-1})(x_j - x_{j-1}) = -f(x_{j-1}), \qquad j = 1, 2, \cdots
\end{equation}
starting from a given initial guess $x_0$. 

\begin{theorem}
\label{th:Convergence}
If the Jacobian $f_x(\bar{x})$ is non-singular in the solution $\bar{x}$ and Lipschitz-continuous
in a neighbourhood of $\bar{x}$, for all $x_0$ sufficiently close to $x$ the sequence
$\lbrace x_j \rbrace$ of the solutions of \eqref{eq:Newton} converges not less than quadratically
to $\bar{x}$.  
\end{theorem}
\begin{proof}
This is a well-known result, see e.g. \cite{Ortega1970}. 
\end{proof}

\begin{theorem}
\label{th:LinearConvergence}
If Equation \eqref{eq:ImplicitEq} is linear and $f_x$ is non-singular, then NR's algorithm
converges in one step, irrespective of the chosen initial guess $x_0$.
\end{theorem}
\begin{proof}
If Equation \eqref{eq:ImplicitEq} is linear, $f(x) = Jx + b$, where $J = f_x$ is a constant
$m \times m$ matrix and $b = f(0)$. The first iteration of \eqref{eq:Newton} becomes
\begin{equation}
J(x_1-x_0) = -(J x_0 + b)
\end{equation}
whose solution $x_1$ is the solution of $Jx + b = 0$.
\end{proof}

\begin{theorem}
\label{th:MixedConvergence}
If NR's algorithm is initialized with a first guess
\begin{equation}
x_0 = \begin{bmatrix} w_0 \\ z_0 \end{bmatrix},
\end{equation}
the values of the approximated solution $x_j$ at each step $j > 0$ only depend on the guess values $w_0$ of the
variables affecting the Jacobian, regardless of the choice of $z_0$.
\end{theorem}
\begin{proof}
Equation \eqref{eq:JacDependency} implies that the Jacobian matrix $f_x$ can be partitioned as
follows
\begin{equation}
\label{eq:JacobianPartition}
f_x(x) = \begin{bmatrix} f_w(w) & f_z \end{bmatrix} , 
\end{equation}
where $f_z$ is a constant matrix. Therefore, the nonlinear function $f$ can be
rewritten as
\begin{equation}
\label{eq:nonLinearFunc2}
f \left(\begin{bmatrix} w \\ z \end{bmatrix} \right) = g(w) + f_z z.
\end{equation}
The first iteration of NR's algorithm \eqref{eq:Newton} yields
\begin{equation}
\label{eq:NewtonStep1}
 \begin{bmatrix} f_w(w_0) &  f_z \end{bmatrix} \begin{bmatrix} w_1 - w_0 \\ z_1 - z_0 \end{bmatrix} =-f\left(\begin{bmatrix} w_0 \\ z_0 \end{bmatrix} \right), 
\end{equation}
which can be expanded into
\begin{equation}
\label{eq:NewtonStep2}
 f_w(w_0) (w_1- w_0) +  f_z (z_1 - z_0) = -g(w_0) - f_z z_0.
\end{equation}
Now the two terms $-f_z z_0$ on the left and right-hand side cancel out, yielding
\begin{equation}
 f_w(w_0) (w_1- w_0) +  f_z z_1  = -g(w_0),
\end{equation}
whose solution
\begin{equation}
\label{eq:NewtonStep4}
x_1 = \begin{bmatrix} w_1 \\ z_1 \end{bmatrix},
\end{equation}
which is the result of the first iteration, does not depend on the initial guess $z_0$. Hence, the values of the subsequent iterations $x_2, x_3, \ldots$ also do not depend on $z_0$.
\end{proof}

\begin{theorem}
\label{th:LinearResiduals}
The residuals of the linear equations in system \eqref{eq:ImplicitEq} after the first iteration of NR's algorithm are zero, i.e, $l(x_1) = 0$, regardless of the initial guess values $x_0$.
\end{theorem}
\begin{proof}
The Jacobian matrix of function \eqref{eq:ImplicitEq} can be partitioned as follows:
\begin{equation}
\label{eq:JacobianLNLPartition}
f_x(x) = \begin{bmatrix} n_w(w) & n_z \\ l_w & l_z \end{bmatrix},
\end{equation}
where the Jacobians $n_z$, $l_w$, and $l_z$ are constant matrices.
The linear equations residuals $l(x)$ can then be formulated as 
\begin{equation}
\label{eq:LinearResiduals}
l(x) = l_w w + l_z z + l(0)
\end{equation}
The first iteration of NR's algorithm \eqref{eq:Newton} reads
\begin{equation}
\label{eq:NewtonStep3}
\begin{bmatrix} n_w(w_0) & n_z \\ l_w & l_z \end{bmatrix} \begin{bmatrix} w_1 - w_0 \\ z_1 - z_0 \end{bmatrix} =-\begin{bmatrix} n(x_0) \\ l(x_0) \end{bmatrix}.
\end{equation}
By expanding and rearranging the last rows of Equation \eqref{eq:NewtonStep3}, and by taking into account Equation \eqref{eq:LinearResiduals} with $x = x_0$, one gets
\begin{equation}
l_w w_1 + l_z z_1 = l_w w_0 + l_z z_0 - l(x_0) = -l(0).
\end{equation}
Hence, 
\begin{equation}
l(x_1) = l_w w_1 + l_z z_1 + l(0) = -l(0) + l(0) = 0.
\end{equation}
\end{proof}

\begin{definition}
Consider the NR algorithm \eqref{eq:Newton}. Assume that the  function $f(x)$ is three times continuously differentiable in an open neighbourhood $\mathfrak{D}$ containing the initial guess $x_0$ and the result of the first iteration $x_1$. Denote the $i$-th component of function $f(x)$ as $f^i(x)$, its Jacobian matrix with respect to $x$ as $f^i_x(x)$, and its Hessian matrix as $f^i_{xx}(x)$. By means of a Taylor series expansion, one can write
\begin{equation}
\label{eq:TaylorExpansion}
f^i(x_1) = f^i(x_0) + f_x^i(x_0)(x_1 - x_0) + \frac{1}{2}(x_1 - x_0)' f^i_{xx}(x_0)(x_1 - x_0) + h^i(x_1,x_0),
\end{equation}
which implicitly defines the higher-order residual functions $h^i(\cdot,\cdot)$. 
\end{definition}

\begin{definition}
With reference to NR's iteration \eqref{eq:Newton}, define the nonlinear residual at iteration point $x_{k-1}$ as
\begin{equation}
\label{eq:NonlinearResidual}
r(x_{k-1}) = f(x_{k-1}) + f_z(z_{k} - z_{k-1})
\end{equation}
\end{definition}

\begin{definition}
Define the coefficients $\alpha_i > 0$, $i = 1, \cdots, m$, such that
\begin{equation}
\label{eq:ConditionOnH}
|h^i(x_1, x_0)| = \alpha_i \norm{r(x_{0})}_\infty,
\end{equation}
and let 
\begin{equation}
\label{eq:MaxAlpha}
\alpha = max(\alpha_i).
\end{equation}
\end{definition}

\begin{definition}
Define the \emph{curvature factor} $\Gamma_{ijk}$ of the $i$-th nonlinear equation with respect to variables $w_j, w_k$ after the first iteration as
\begin{equation}
\label{eq:GammaDefinition}
\Gamma_{ijk} = \abs{\frac{1}{2} \frac{\partial^2 g^i(w_0)}{\partial w_j \partial w_k} \frac{ (w_{1,k} - w_{0,k})(w_{1,j} - w_{0,j})}{\norm{r(x_0)}_\infty}}, \qquad
\begin{matrix} i=1,\ldots,p, \\ j=1,\ldots,q, \\ k=1,\ldots,q. \end{matrix}
\end{equation}
\end{definition}

\begin{theorem}
\label{th:SufficientCondition}
Given a constant $\beta > 0$, a sufficient condition for the property 
\begin{equation}
\label{eq:SufficientProperty}
\norm{f(x_1)}_\infty \leq (\alpha + \beta) \norm{r(x_0)}_\infty 
\end{equation}
to hold is that 
\begin{equation}
\label{eq:SufficientCondition}
\sum_{jk} \Gamma_{ijk} \leq \beta \qquad \forall i = 1, \cdots, p.
\end{equation}
\begin{proof}
Equation \eqref{eq:Newton} for the first iteration reads
\begin{equation}
\label{eq:NewtonFirstStep}
f_x(x_0)(x_{1} - x_0) = -f(x_0)
\end{equation}
Computing $f^i(x_1)$ with Equation \eqref{eq:TaylorExpansion} and plugging in Equation \eqref{eq:NewtonFirstStep}, one obtains
\begin{equation}
\label{eq:IterationResidual}
f^i(x_1) =  \frac{1}{2}\sum_{jk}{\frac{\partial^2 f^i(x_0)}{\partial x_j \partial x_k}(x_{1,j} - x_{0,j})(x_{1,k} - x_{0,k})}  + h^i(x_1, x_0).
\end{equation}
Recalling \eqref{eq:JacobianPartition}, \eqref{eq:ConditionOnH}, \eqref{eq:MaxAlpha},  \eqref{eq:GammaDefinition}, the following chain of inequalities holds:
\begin{align}
|f^i(x_1)| & \leq \abs{\frac{1}{2}\sum_{jk}{\frac{\partial^2 f^i(x_0)}{\partial x_j \partial x_k}(x_{1,j} - x_{0,j}) (x_{1,k} - x_{0,k})} + h^i(x_1,x_0)}\\
& \leq \abs{\frac{1}{2}\sum_{jk}{\frac{\partial^2 g^i(w_0)}{\partial w_j \partial w_k}(w_{1,j} - w_{0,j}) (w_{1,k} - w_{0,k})} + h^i(x_1, x_0)}\\
& \leq \sum_{jk}\abs{\frac{1}{2}{\frac{\partial^2 g^i(w_0)}{\partial w_j \partial w_k}(w_{1,j} - w_{0,j}) (w_{1,k} - w_{0,k})}} + \abs{h^i(x_1, x_0)}\\
& \leq \sum_{jk}\Gamma_{ijk}\norm{r(x_0)}_\infty + \alpha_i \norm{r(x_0)}_\infty \\
& \leq \beta \norm{r(x_0)}_\infty + \alpha \norm{r(x_0)}_\infty \\
& \leq (\alpha+\beta)\norm{r(x_0)}_\infty.
\end{align}

By recalling the definition of the $\infty$-norm, it follows immediately that 

\begin{equation}
\norm{f(x_1)}_\infty = \max_i{\abs{f^i(x_1)}} \leq (\alpha+\beta)\norm{r(x_0)}_\infty.
\end{equation}

\end{proof}
\end{theorem}

\begin{remark}
$r^i(\bar{x}) = r^i(\bar{w}) = f^i(\bar{x}) = 0$
\end{remark}

\begin{remark}
If the Jacobian $J(x)$ is non-singular, then $(w_{1,k} - w_{0,k}) = 0 \quad \forall k \Leftrightarrow f(x_1) = 0$.
\end{remark}

\begin{remark}
The coefficients $\alpha_i$ of the nonlinear equations, showing up in Theorem \ref{th:SufficientCondition}, can be computed from Equations \eqref{eq:Newton},\eqref{eq:JacobianPartition},\eqref{eq:TaylorExpansion}, and \eqref{eq:ConditionOnH}, yielding
\begin{align}
\alpha_i & = \frac{\abs{f^i(x_1) - \frac{1}{2}(x_1 - x_0)' f^i_{xx}(x_0)(x_1 - x_0)}}{\norm{r(x_0)}_\infty} \\
\label{eq:Alpha_i}
& = \frac{\abs{f^i(x_1) - \frac{1}{2}(w_1 - w_0)' f^i_{ww}(w_0)(w_1 - w_0)}}{\norm{r(x_0)}_\infty}
\end{align}
\end{remark}

\begin{remark}
$\alpha_i = 0$ for $i = p+1, \cdots, m$, since the Taylor expansion of linear equations obviously lacks all terms of order greater than one. Furthermore, $\alpha_i = 0$ for all quadratic equations.
\end{remark}

\begin{remark}
Assuming the system of equations \eqref{eq:ImplicitEq} comes from a physical modelling problem, both the unknowns $w$ and the residuals $f^i(x)$ are in general dimensional quantities. The factors $\alpha_i$ as defined in \eqref{eq:ConditionOnH} are obviously non-dimensional and invariant with respect to the choice of all units of the problem. The curvature factors defined in \eqref{eq:GammaDefinition} are invariant with respect to the choice of unit of the unknowns $w$, which appear both at the numerator and at the denominator. However, they are not invariant with respect to the choice of units of the equation residuals $f^i(x)$. Also, computing the $\infty$-norm of the residual $f(x)$ is conceptually questionable, as it involves finding the maximum of quantities with different dimensions (what is the maximum between one meter and ten watts?) and turns out to be more sensitive to those components that have large numerical values because of the choice of a small physical unit. It is therefore recommended to scale the residual functions with appropriately sized quantities, so that all the residuals $f^i(x)$ are non-dimensional and have about the same orders of magnitude, see, e.g., \cite{CasellaBraunEOOLT2017} for a method to accomplish this.
\end{remark}

\begin{theorem}
\label{th:Invariance}
Given a problem \eqref{eq:ImplicitEq} and an initial guess $w_0$ for the nonlinear variables, the values of $r^i(x_0)$, $\alpha_i$, and $\Gamma_{ijk}$ are invariant with respect to the choice of the initial guess $z_0$ for the linear variables. 
\end{theorem}
\begin{proof}
The first NR iteration can be expanded as
\begin{equation}
f(x_0) + f_w(w_0)(w_1 - w_0) + f_z(z_1 - z_0) = 0,
\end{equation}
which implies
\begin{equation}
\label{eq:InvarianceOfR}
r(x_0) = f(x_0) + f_z(z_1 - z_0) = - f_w(w_0)(w_1 - w_0).
\end{equation}
Hence, from Equation \eqref{eq:Alpha_i}, $\alpha_i$ can be computed as
\begin{equation}
\alpha_i = \frac{\abs{f^i(x_1) - \frac{1}{2}(w_1 - w_0)' f^i_{ww}(w_0)(w_1 - w_0)}}{\norm{f_w(w_0)(w_1 - w_0)}_\infty}.
\end{equation}
According to Theorem \ref{th:MixedConvergence}, $x_1$ and $w_1$ do not depend on $z_0$, therefore also $\alpha_i$ doesn't. Considering Equation \eqref{eq:InvarianceOfR} again, it is apparent from Equation \eqref{eq:GammaDefinition} that $\Gamma_{ijk}$ also does not depend on $z_0$.
\end{proof}

\begin{theorem}
\label{th:Sensitivity}
The sensitivity of the solution $x_1$ after the first NR iteration, with respect to changes in the initial guess $x_0$, can be computed as:
\begin{equation}
\label{eq:Sensitivity}
\frac{\partial x_1}{\partial x_0} = \Sigma,
\end{equation}
where
\begin{align}
\label{eq:HtildeStart}
& H_i =  (w_1-w_0)' f^i_{ww}(w_0) \\
& H = \begin{bmatrix} H_1 \\ H_2 \\ \cdots \\ H_p \end{bmatrix} \\
\label{eq:DefinitionOfS}
& \Sigma =  -\left[{f_x(w_0)} \right] ^{-1} \begin{bmatrix} H_{p\times q} & 0_{p
\times (m-q)} \\ 0_{(m-p) \times q} & 0_{(m-p) \times (m-q)} \end{bmatrix}.
\end{align}
\end{theorem}

\begin{proof}
The first NR iteration, Equation \eqref{eq:Newton} with $j = 1$, reads
\begin{equation}
f_x(x_0)(x_1 - x_0) = -f(x_0)
\end{equation}
By differentiating the $i$-th row with respect to $x_0$, one obtains the following $1 \times m$ matrix equation
\begin{equation}
 (x_1 - x_0)' f^i_{xx}(x_0) + f^i_x(x_0)\frac{\partial(x_1 - x_0)}{\partial x_0} = -f^i_x(x_0),
\end{equation}
where $f^i_x$ and $f^i_{xx}$ are the Jacobian and Hessian matrices of the i-th equation residual function. By stacking the $m$ row vectors corresponding to each equation in \eqref{eq:ImplicitEq} and recalling that all the derivatives of function $f(x)$ only depend on $w$, one obtains the following matrix equation with $m \times m$ terms
\begin{equation}
\label{eq:Stacked}
\tilde{H} + f_x(w_0) \frac{\partial(x_1 - x_0)}{\partial x_0} = -f_x(w_0),
\end{equation}
where 
\begin{align}
\label{eq:Htildei}
& \tilde{H}_i =  (x_1 - x_0)' f^i_{xx}(w_0) \\
& \tilde{H} = \begin{bmatrix} \tilde{H}_1 \\ \tilde{H}_2 \\ \cdots \\ \tilde{H}_m \end{bmatrix}
\end{align}
which can be solved for the sensitivity matrix by left-multiplying each term in equation \eqref{eq:Stacked} by the inverse Jacobian $\left[{f_x(w_0)} \right] ^{-1}$, yielding
\begin{equation}
\label{eq:SensitivityProof1}
\frac{\partial}{\partial x_0}(x_1 - x_0) = - \left[{f_x(w_0)} \right] ^{-1} \tilde{H} - I_{m \times m}.
\end{equation}
Considering that
\begin{equation}
\frac{\partial(x_1 - x_0)}{\partial x_0} = \frac{\partial x_1}{\partial x_0} - I_{m \times m},
\end{equation}
equation \eqref{eq:SensitivityProof1} can be reduced to 
\begin{equation}
\frac{\partial x_1}{\partial x_0} = - \left[{f_x(w_0)} \right] ^{-1} \tilde{H}.
\end{equation}
Since the first derivatives of $f(x)$ only depend on the first $p$ elements of
vector $x$ (i.e., the $w$ vector), the last $m-q$ rows and columns of
$f^i_{xx}$ are zero. Hence, the last $m-q$ columns of each $\tilde{H}_i$ are
zero, and so are the last $m-q$ columns of the stacked matrix $\tilde{H}$.
Furthermore, when computing the matrix product inside \eqref{eq:Htildei}, the
last $m-q$ terms of the $x_1 - x_0$ vector, i.e., $z_1 - z_0$, always get
multiplied by zero second derivatives, so they can be skipped. Finally, since
the last $m-p$ equation residuals are linear, their Hessians are zero, so the
last $m-p$ row vectors  $\tilde{H}_i$ are also zero.

Hence, it is possible to compute matrix $\tilde{H}$ more efficiently by skipping all those elements that do not contribute to the final result, yielding
\begin{align}
& H_i = (w_1-w_0)' f^i_{ww}(w_0), \qquad i = 1,\cdots,p \\
& H = \begin{bmatrix} H_1 \\ H_2 \\ \cdots \\ H_p \end{bmatrix} \\
& \tilde{H} =  \begin{bmatrix} H_{p\times q} & 0_{p \times (m-q)} \\ 0_{(m-p) \times q} & 0_{(m-p) \times (m-q)} \end{bmatrix}
\end{align}
\end{proof}

\begin{remark}
When computing the sensitivity \eqref{eq:Sensitivity} at the solution $x_0 = \bar{x}$, since $w_1 = w_0 = \bar{w}$, it follows that $H = 0$, so the sensitivity $\Sigma$ turns out to be zero. 

This means that if an initial guess equal to the solution plus an
infinitesimally small perturbation $x_0 = \bar{x} + \delta x$ is chosen, the
solution $x_1$ after the first NR iteration is not affected at all. This is
consistent with the fact that $f(x)$ can be approximated as a linear function
in a small neighbourhood of the solution $\bar{x}$, so that Theorem \ref{th:LinearConvergence}
guarantees that the first NR iteration converges to the solution $\bar{x}$ in
just one iteration, irrespective of the initial guess.

If the initial guess $x_0$ is close enough to the solution $\bar{x}$ that the function $f(x)$ is still approximately linear in a neighbourhood containing $x_0$ and $\bar{x}$, then the same behaviour is preserved, i.e., the result after the first iteration is insensitive to small changes of the initial guess, so $\Sigma \approx 0$.

As $x_0$ is chosen farther away from the solution $\bar{x}$, nonlinear effects kick in, accounted for by matrix $\Sigma$, which can be then considered an indicator of how far the initial guess is from the sweet spot of NR convergence. 

In particular, a value $\abs{\sigma_{jj}} \ll 1$ means that the effect of applying a certain perturbation to $w_{0,j}$ on $w_{1,j}$ will be less than the perturbation itself, meaning that the nonlinear effects are moderate, while $\abs{\sigma_{jj}} \gg 1$ means that the effect of a certain perturbation is amplified after the first iteration, which is a sure sign of large nonlinear effects on the NR iteration pertaining to that variable.

Indeed, according to Equations \eqref{eq:HtildeStart}-\eqref{eq:DefinitionOfS}, matrix $\Sigma$ becomes larger as the initial guesses of the nonlinear variables $w_0$ get farther away from $\bar{w}$, increasing $(w_1 - w_0)$; this also depends on how large the corresponding second derivatives  in the Hessians $f^i_{ww}$ are. On the other hand, the fact that the initial value of the linear variables $z_0$ gets farther from $\bar{z}$ is completely irrelevant, since $(z_1 - z_0)$ does not affect the value of $\Sigma$ at all. This is also consistent with Theorem \ref{th:MixedConvergence}. 
\end{remark}

\begin{remark}
\label{rem:SensitivityStructure}
The last $m-q$ columns of matrix $\Sigma$ are zero; this means that the sensitivity of the increment $x_1 - x_0$ with respect to $z_0$ is nil; this is also consistent with Theorem \ref{th:MixedConvergence}. However, the sensitivity of $z_1 - z_0$ with respect to $w_0$ can in general be non-zero. Hence, matrix $\Sigma$ has the following structure:
\begin{equation}
\label{def:Sigma}
\Sigma = \begin{bmatrix} \Sigma_{ww} & 0_{q \times (m-q)} \\ \Sigma_{zw} & 0_{(m-q)\times(m-q)} \end{bmatrix}
\end{equation}
\end{remark}

\begin{remark}
\label{rem:SensitivityScaling}
Assuming that the system of equations \eqref{eq:ImplicitEq} comes from a physical modelling problem, the dimension of a generic element $\sigma_{j,k}$ of matrix $\Sigma$ is the dimension of $x_j$ divided by the dimension of $x_k$. Hence, the diagonal elements $\sigma_{jj}$ of $\Sigma$ are non-dimensional and scale invariant, while the off-diagonal elements in general are not, so that off-diagonal terms $\sigma_{jk}$ are in general not invariant with respect to a change of units in the formulation of the physical problem. This also means that the fact that one such element $\sigma_{jk}$ is much smaller or much greater than one doesn't have any particular meaning, since its actual value depends on the choice of units of the problem, which is completely arbitrary, and is not invariant with respect to the scaling of the problem.
\end{remark}

%


\section{Discussion}
\label{sec:Discussion}
The theorems stated in the previous section can be used to formulate four criteria for the selection of the initial guess values $x_0$ for NR's algorithm.

The well-know Theorem \ref{th:Convergence} implies that if the initial guess $x_0$ is close enough to the sought solution $\bar{x}$, NR's algorithm converges quickly to the exact solution; however, it does not provide any indication on how close the initial guesses $x_0$ must be to the solution $\bar{x}$ for this outcome to take place. 

Theorem \ref{th:LinearConvergence} indicates that in case function $f(x)$ is fully linear, NR's algorithm always converges in one step, no matter what the initial guess $x_0$ is. This is an interesting limit case, but it hardly has any practical importance, since NR algorithms are normally employed to solve systems that include nonlinear equations.

Theorem \ref{th:MixedConvergence}, instead, is of much greater practical importance when dealing with systems of mixed linear and nonlinear equations, a case often encountered in applications, as it states that the value of $x_1$ after one iteration will be the same regardless of the initial guess of the linear variables $z_0$. Hence, taking care of providing initial guess values for $z$ is a complete waste of time; one should rather invest time and effort in providing good initial guesses $w_0$ for the nonlinear variables $w$ that actually influence the Jacobian. 

This consideration is valid assuming that a direct method, e.g. LU decomposition) is used to solve the linear system; in case an iterative method is used, the initial guess of the linear part could play a role in determining the number of iterations of the linear solver, and thus the performance of NR's algorithm. However, the typical size of the problems addressed by this paper does not usually exceed the order of magnitude of tens of thousands; besides, the structure of the problems is usually characterized by a high degree of sparsity, with only a handful of variables showing up in each equation. Direct solvers are normally employed in these cases, possibly using sparse algorithms such as KLU \cite{Davis2010} if $p$ exceeds a few tens of equations, with satisfactory performance, see e.g. \cite{BraunEtAlModelica2017}, making the statement above valid in practical applications.

Based on these considerations, the following first Criterion can then be formulated:

\begin{criterion}
\label{cr:InitialGuesses}
When choosing the initial guesses $x_0$ for NR's algorithm, provide good initial guesses for the variables $w$ that influence the Jacobian. The other variables $z$, that only appear linearly in the system of equations, can be given a trivial initial guess $z_0 = 0$, without affecting the convergence of NR's algorithm. 
\end{criterion}

If the initial guess $x_0$ is close enough to the solution $\bar{x}$, the convergence of the sequence $x_p$ is fast and the norm of the residual $\norm{f(x)}_\infty$ becomes much smaller than the previous one at each iteration $p$, already starting from the first one. A good initial guess $x_0$ could then be identified as fulfilling the following property:
\begin{equation}
\label{eq:ResidualReduction}
\norm{f(x_1)}_\infty \ll \norm{f(x_{0})}_\infty.
\end{equation}
However, if one wants to exploit the invariance of the NR iterations with respect to the initial guess of the linear variables $z_0$  (Theorem \ref{th:Invariance})  and thus apply Criterion \ref{cr:InitialGuesses}, the trivial initial guess $z_0 = 0$ may be quite far from the solution $\bar{z}$, possibly causing $\norm{f(x_0)}_\infty$ to become very large. In this case, a large reduction of the norm of the residual after the first iteration may simply be the effect of the residuals of some linear equations becoming zero after the first iteration because of Theorem \ref{th:LinearResiduals}, rather than being an indication that the initial values of the nonlinear variables are close to the solution.

A better indication of closeness to convergence can then be obtained by deducting from $f(x_0)$ the effect of the linear variables increment $(z_1 - z_0)$, hence using the nonlinear residuals $r(x_0)$ as defined in Equation \eqref{eq:NonlinearResidual} in place of $f(x_0)$. The following condition is then sought:
\begin{equation}
\label{eq:SufficientProperty2}
\norm{f(x_1)}_\infty \ll \norm{r(x_0)}_\infty,
\end{equation}
which has the nice property of not depending on the choice of $z_0$ thanks to Theorem \ref{th:Invariance}, and thus to be fully consistent with Criterion \ref{cr:InitialGuesses}.

Note that condition \eqref{eq:SufficientProperty2} still only has heuristic value, as it is possible to build counter-examples where this property holds, but then subsequent iterations do not converge to any solution. However, in most practical cases, if one makes an effort to choose a good initial guess and condition  \eqref{eq:SufficientProperty2} is fulfilled, it is quite unlikely that convergence is not eventually achieved.

Conversely, if despite all efforts, convergence to the solution is \emph{not} achieved and condition \eqref{eq:SufficientProperty2} is \emph{not} satisfied, the likely explanation for the convergence failure is that at least one of the initial guesses $w_0$ is not close enough to the solution, and as a consequence the curvature of the hyper-surfaces $y = f(x)$ causes the increment $(w_1 - w_0)$ to go astray, blowing up the nonlinear residual after the first iteration. In this case, one can formulate heuristic criteria based on  Theorems \ref{th:SufficientCondition} and \ref{th:Sensitivity} to understand which components of the vector of initial guesses $w_0$ are most likely responsible for this, and should then be improved to eventually achieve convergence.

The sufficient condition of Theorem \ref{th:SufficientCondition} to obtain property \eqref{eq:SufficientProperty2} requires 
\begin{equation}
\alpha + \beta \ll 1
\end{equation}
to hold. If $\alpha_i > 1$ for some $i$, then the sufficient condition is
violated. In this case, the problem is that across the first iteration of NR's
algorithm, the higher-order residual $h^i(w_1, w_0)$ of the Taylor expansion \eqref{eq:TaylorExpansion} plays a major role, which obviously contradicts the requirement that all iterations
should take place in a neighbourhood of the solution, where the functions
$f^i(x)$ are well approximated by linear ones, ensuring fast convergence.
Hence, the initial values of the nonlinear variables appearing in $f^i(x)$ are
not good enough, and should be improved. Unfortunately, it is not possible in
this case to discriminate among the role played by each individual nonlinear
variable, in case more than one is involved in the $i$-th equation.

In case $\alpha \ll 1$, then the sufficient condition is satisfied if also $\beta \ll 1$, which can only be achieved if $\Gamma_{ijk} \ll 1$,  $\forall i,j,k$. If a certain $\Gamma_{ijk} > 1$, this causes the violation of the sufficient condition, which may be due to the fact that poor initial guesses of the $j$-th and $k$-th nonlinear variables in vector $w_0$ cause the first iteration of NR's algorithm to span an interval where the curvature of the corresponding hyper-surfaces in $y = f^i(x)$ is large enough to potentially cause convergence problems.

Note that a violation of the sufficient (but not necessary!) condition of Theorem \ref{th:SufficientCondition} does not necessarily mean that the residual norm $\norm{f(x_s)}_\infty$ will not get smaller with increasing iteration number $s$, nor that NR's algorithm will not eventually converge. However, in case of convergence failure of NR's algorithm, Theorem \ref{th:SufficientCondition} can provide useful indications about which components of initial guesses in $w_0$ may be responsible for the failed convergence.

In case the initial guess $w_0$ is sufficiently far from the solution $\bar{w}$, it may happen that vector $x_1$ of the unknowns after the first iteration does not belong to the domain of definition of the residual function $f(x)$, preventing the computation of $f(x_1)$ and thus the computation of the $\alpha_i$ factors. In order to still obtain useful information about the variables potentially causing the convergence failure because of high-order terms in the Taylor expansion, one can compute a \emph{damped} first NR iteration $x_1^{*}$ such that
\begin{equation}
\label{eq:DampingNewton}
x_1^{*} - x_0 = \lambda (x_1 - x_0), \quad 0 < \lambda \leq 1
\end{equation}
or, equivalently
\begin{equation}
\label{eq:DampedNewton}
f_x(x_0)(x_1^{*} - x_0) = -\lambda f(x_0).
\end{equation}

By taking a small enough value of $\lambda$, one can get $x_1^{*}$ arbitrarily close to the initial guess $x_0$, hence within the domain of definition of $f(x)$, assuming that $x_0$ is an interior point of that domain. One can then exploit the Taylor expansion
\begin{equation}
\label{eq:TaylorExpansionDamped}
f^i(x_1^{*}) = f^i(x_0) + f_x^i(x_0)(x_1^{*} - x_0) + \frac{1}{2}(x_1^{*} - x_0)' f^i_{xx}(x_0)(x_1^{*} - x_0) + h^i(x_1^{*}, x_0)
\end{equation}
to compute $h^i(x_1^{*}, x_0)$. Is it then possible to re-define $\alpha_i$ as
\begin{equation}
\label{eq:NewAlphaDefinition}
|h^i(x_1^{*}, x_0)| = \alpha_i \lambda^3 \norm{r(w_0)}_\infty,
\end{equation}
where $\lambda^3$ accounts for the fact that the term $h^i(x_1^{*}, x_0)$ shrinks as $\lambda^3$ asymptotically as $\lambda \rightarrow 0$, thus making definition \eqref{eq:NewAlphaDefinition} 
asymptotically invariant as $\lambda \rightarrow 0$. By combining the previous three equations and by taking into account that only the $w$ unknowns are relevant for the quadratic term, one can compute $\alpha_i$ as:
\begin{equation}
\label{eq:Alpha_i_Modified}
\alpha_i = \frac{\abs{f^i(x_1^{*}) - (1 - \lambda)f(x_0) - \frac{1}{2}(w_1^{*} - w_0)' f^i_{ww}(w_0)(w_1^{*} - w_0)}}{\lambda^3 \norm{r(w_0)}_\infty} 
\end{equation}
One should then reduce $\lambda$ until $x_1^{*}$ is close enough to $x_0$ to allow computing $f(x_1^{*})$, then use Equation \eqref{eq:Alpha_i_Modified} to compute the $\alpha_i$ factors.

Summing up, the idea is that the larger values of $\alpha_i$ and/or $\Gamma_{ijk}$ point to the initial guesses which are possibly the cause of the convergence failure. If those initial guesses are found not to correspond to a reliable prior estimate, e.g. because of some programming error or oversight, they can be fixed according to that. Otherwise, one can use the signs of the increments $(w_1 - w_0)$ to get some indication whether the initial guesses should be increased or reduced, though those increments do not give any reliable information about the required magnitude of such a change. 

All of these results can be summarized in the following Criterion.

\begin{criterion}
\label{cr:AlphaGamma}
In case of failure of NR's algorithm to converge to the desired solution $\bar{x}$ starting from the initial guess $w_0$, focus the attention on the initial guess of the variables $w$ that appear in nonlinear equations with $\alpha_i > 1$, as well as on the initial guess of those variables $w$ corresponding to the indices $j$ and $k$ of the curvature factors $\Gamma_{ijk} > 1$, computed after the first iteration of NR's algorithm. Larger values of $\alpha_i$ and $\Gamma_{ijk}$ are likely to correspond to more critical values. The initial guesses may be improved by increasing or decreasing them according to the sign of their increment after the first NR iteration.
\end{criterion}

This criterion implicitly assumes that if the initial guess of the $j$-th nonlinear variable $w_{0,j}$ is far enough from the solution $\bar{w}_{j}$, this will affect the corresponding $\alpha_i$ and $\Gamma_{ijk}$ indicators \emph{and only them}, so that they can be used backwards to pinpoint the critical initial guesses. In principle, this is not necessarily the case: it is possible that a significant error on $w_{0,j}$ has an influence on $w_{1,k}, k \neq j$, and thus on the $\alpha_i$ and $\Gamma_{ijk}$ indicators pertaining to the $k$-th nonlinear variable, leading to the potentially incorrect diagnosis that the initial guess $w_{0,k}$ is wrong too and needs to be improved. In other words, an error on the initial guess of the $j$-th nonlinear variable could spill over to the indicators of other variables, leading to false positives of Criterion \ref{cr:AlphaGamma}.

One way to spot this potential spill-over effect and be alerted about possible false positives of Criterion \ref{cr:AlphaGamma} could be to look at the off-diagonal elements of sensitivity matrix $\Sigma$ introduced in \eqref{def:Sigma}, to check if and by how much an error on the initial guess of the $j$-th nonlinear variable can have an influence on the $k$-th nonlinear variable after the first iteration. Unfortunately, as already noted in Remark \ref{rem:SensitivityScaling}, the magnitude of such off-diagonal elements depends crucially on the (arbitrary) scaling of the corresponding variables. The authors tried several scaling methods, but none of them brought consistently reliable indications in all the test cases, leading to conclude that the off-diagonal elements of $\Sigma$ serve no useful purpose in the context of this paper.

On the other hand, the diagonal elements of matrix $\Sigma$ can be also be used to provide information about the initial guesses $w_{0,j}$ which are most likely the cause of the NR convergence failure. As long as the initial guess $w_0$ is close enough to the solution $\bar{w}$ that the function $f(x)$ is approximately linear in a neighbourhood containing both, then Theorems \ref{th:Convergence} and \ref{th:LinearConvergence} indicate that $x_1 \approx \bar{x}$, irrespective of the actual value of $x_0$; hence, the sensitivity of $x_{1,j}$ to $x_{0,j}$ will be very small, i.e., $\abs{\sigma_{jj}} \ll 1 \;\; \forall j$. Conversely, an element $\abs{\sigma_{jj}} > 1$ indicates that a small change on the initial guess $w_{0,j}$ has an effect of larger magnitude on the result of the first iteration $w_{1,j}$, which is incompatible with $w_{0,k}$ being close enough to the solution to be in the sweet spot of superlinear convergence. Hence, a large value of $\sigma_{jj}$ suggests to check the initial guess $w_{0,j}$, possibly trying to improve it based on the sign of the corresponding variable increment after the first NR iteration.

From this point of view, the elements $\sigma_{jj}$ provide second-order information which is somewhat related, but at the same time complementary, to the second-order information provided by the $\Gamma_{ijk}$ indicators.

On the positive side, this information is based on the joint effects of all the unknowns on the convergence process, so it doesn't suffer from the spillover effect that may cause false positives if one only looks at $\alpha_i$ and $\Gamma_{ijk}$. It is also invariant with respect to the scaling of the problem.

These considerations can be summarized in the following criterion:

\begin{criterion}
\label{cr:Sigma}
In case of failure of NR's algorithm to converge to the desired solution $\bar{x}$ starting from the initial guess $w_0$, focus the attention on the initial guess of the variables $w_j$ such that $\abs{\sigma_{jj}} > 1$. Larger values are likely to indicate a stronger effect on the convergence failure. The initial guesses can be improved based on prior knowledge, or according to the sign of their increment after the first NR iteration.
\end{criterion}

On the other hand, there are several possible shortcomings in the use of Criterion \ref{cr:Sigma}. First of all, matrix $\Sigma$ relies exclusively on local first- and second-order information about function $f(x)$ and has no provision to take into account higher-order effects, contrary to the $\alpha_i$ indicators of Criterion \ref{cr:AlphaGamma}. 

Secondly, while $\Gamma_{ijk}$ is sparse, because there are usually only a handful of pairs of unknowns ($w_j$, $w_k$) showing up in the $i$-th nonlinear equation in most physical models, the $\Sigma$ matrix is not, so that its computation can become a bit cumbersome as the size of the problem grows, even though the most computationally intensive step is the factorization of the Jacobian matrix, which needs to be carried out anyway for the NR algorithm iteration. 

Last, but not least, $\Sigma$ provides information about potentially wrong initial guesses, but gives no indication on which equations are possibly causing trouble. This information could be relevant for model developers, who may consider using equivalent formulations of those equations, that lead to the same solution but are less critical from the point of view of nonlinear behaviour.

One important consideration is due at this point. The theorems proven in Section \ref{sec:Method} lead to the conclusion that if $\alpha_i \ll 1$, $\Gamma_{ijk} \ll 1$, and $\abs{\sigma_{jj}} \ll 1$ $\forall i, j, k$, then the initial guess is already most likely in the sweet spot of superlinear convergence of NR's method. However, those conditions are by no means necessary to achieve convergence, hence one should not try to fulfil all of them, because that will most likely end up in an unnecessary waste of time. Furthermore, there are no rigorous arguments to establish what are the upper limits for $\alpha_i$, $\Gamma_{ijk}$, and $\abs{\sigma_{jj}}$, above which one should focus on the corresponding variable(s); in Criteria \ref{cr:AlphaGamma} and \ref{cr:Sigma} a limit of one was suggested, but that limit value is in fact somewhat arbitrary.

One could then formulate the following heuristic argument, that those indicators should be ranked in descending order, and one should first focus on the largest one(s), try to fix the corresponding variable(s), and then iterate until the initial guess turns out to be good enough to achieve convergence, even if all the indicators are not well below one. This argument is particularly reasonable under the assumption that some effort has been made to provide good enough initial guesses, so only a handful of variables may still be problematic from the point of view of convergence. Although there is no rigorous way to back this statement, the examples shown in Sect. \ref{sec:Examples} demonstrate that this strategy can be quite effective. Based on this argument, one can then formulate the following Criterion

\begin{criterion}
\label{cr:Ranking}
In case of failure of NR's algorithm to converge to the desired solution $\bar{x}$ starting from the initial guess $w_0$, $z_0 = 0$, compute the $\alpha_i$, $\Gamma_{ijk}$, and $\abs{\sigma_{jj}}$ indicators and rank them in descending order. Then, identify the variable(s) that correspond to the indicators with the largest value and try to improve their initial guess based on prior knowledge or by increasing or decreasing them based on the sign of their increment. If NR's method does not converge to the desired solution with the improved initial guess, compute and rank the indicators one more time, and repeat until convergence is eventually achieved.
\end{criterion}

\section{Example cases}
\label{sec:Examples}
In this section, the practical usefulness and feasibility of the criteria and algorithm presented in the previous section are demonstrated in three example cases, a small thermo-fluid system, an electrical DC circuit and a large AC electrical power system.

\subsection{Thermo-hydraulic system example}
Consider a system comprising a heat exchanger, that absorbs heat from an environment at fixed temperature $T_a$. The working fluid, with specific heat capacity $c$, comes from a source at fixed pressure and temperature $p_s$, $T_s$. 

The fluid first flows through a shut-off valve, which is normally open and thus has a very large flow coefficient $k_p$, corresponding to a small pressure drop, then flows through the heat exchanger, which has a certain pressure drop depending on the coefficient $k_h$, and finally flows through a control valve with flow coefficient $k_v$, that discharges at a fixed pressure $p_d$. The heat duty $Q$ depends on area $A$ of the heat transfer surface, on the specific heat transfer coefficient $\gamma$ that follows a power law depending on the flow rate $w$, and on the difference between the environment temperature and the fluid average temperature.

The system is described by the following equations
\begin{align}
\label{eq:FirstPressureLoss}
0 &= f - k_p \sqrt{p_s - p_i}\\
\label{eq:SecondPressureLoss}
0 &= p_i - p_o - k_h f^2 \\
0 &= f - k_v \sqrt{p_o - p_d} \\
0 &= Q - fc(T_o - T_s) \\
0 &= Q - \gamma A \left( T_a - \frac{T_s+T_o}{2} \right) \\
0 &= \gamma - \gamma_0 {\left(\frac{f}{f_0}\right)}^\nu
\end{align}

The goal of the problem is to find the value of the valve flow coefficient $k_v$ that delivers a certain required heat duty $Q$. Hence, $p_s$, $p_d$, $k_p$, $k_h$, $c$, $f_0$, $\gamma_0$, $\nu$, $T_s$, $T_a$, $Q$, $A$ are known parameters, while $x = w = \begin{bmatrix}f & k_v & T_o & \gamma & p_o & p_i \end{bmatrix}'$. 

Taking $p_s = 2.201$, $p_d = 1$, $k_p = \sqrt{1000}$, $k_h=0.2$, $c = 1$, $f_0 = 1$, $\gamma_0 = 1$, $\nu = 0.8$, $T_s=0$, $T_a = 6$, $Q = 4$, $A = 1$,  the system has an exact solution $f = 1$, $k_v = 1$, $T_o = 4$, $\gamma = 1$, $p_o = 2$, $p_i = 2.2$.


In this example, all variables affect the Jacobian, but they do not do so with the same intensity. The most strongly nonlinear equation is the first one \eqref{eq:FirstPressureLoss}: due to the very large flow coefficient $k_p$, which leads to a very small pressure drop $p_s-p_i$, the first equation residual is much more sensitive to errors in the initial guess of its nonlinear variable $p_i$ than all the other ones. In fact, an error of $1\%$ on the initial guess of $p_i$ can have dramatic consequences on the convergence of NR's algorithm, while all other variables can tolerate initial guesses with errors of $20-30\%$ without substantially hampering the convergence. 

It is expected that the criteria proposed in the previous section allow to get to the same conclusions automatically, without the need of any such expert insight on the mathematical properties of the system. 

Table \ref{tab:HEXConvergence} reports the initial guesses $w_0$, the number of NR iterations, and all the relevant $\alpha_i$, $\Gamma_{ijk}$, and $\Sigma$ indicators, corresponding to different choices of $w_0$.   The NR algorithm is stopped when the absolute value of the largest increment after the last iteration is less than $10^{-12}$.

Since values of the indicators much smaller than unity are not relevant to the analysis, results are displayed with three decimal digits only, to avoid cluttering the presented results with irrelevant detail. Note that the residuals of the second, fourth, and fifth equations are second-degree polynomials, hence $\alpha_2=\alpha_4=\alpha_5=0$ irrespective of the chosen initial guess.

\begin{table}[hbtp]
  \begin{center}
    \caption{Convergence analysis of heat exchanger test}
    \label{tab:HEXConvergence}
    \begin{small}
    \begin{tabular}{c|c|c|c|c|c|c}    
		Var           & \#1       & \#2   & \#3   & \#4     & \#5  & \#6   \\
		\hline
		$f$           & 0.99999  & 0.999    & 0.99   & 0.9     & 0.9    & 3.00  \\
		$k_v$         & 0.99999  & 0.999    & 0.99   & 0.9     & 0.9    & 0.999 \\
		$T_0$         & 3.99996  & 3.996    & 3.96   & 3.6     & 3.6    & 3.996 \\
		$\gamma$      & 0.99999  & 0.999    & 0.99   & 0.9     & 0.9    & 0.999 \\
		$p_o$         & 1.99998  & 1.998    & 1.98   & 1.8     & 1.8    & 1.998 \\
		$p_i$         & 2.19998  & 2.198    & 2.178  & 1.98    & 2.151  & 2.198 \\
		$N_{iter}$    & 3        & 5        & --	    & --      & --     & --    \\
		$\lambda$     & 1.000    & 1.000    & 0.490  & 0.490   & 0.490  & 0.700  \\
		$\alpha_{1}$   & 0.000 & 0.224 & 0.678 & 1.316 & 0.902 & 0.028 \\
		$\alpha_{3}$   & 0.000 & 0.000 & 0.000 & 0.000 & 0.000 & 0.013 \\
		$\alpha_{6}$   & 0.000 & 0.000 & 0.000 & 0.000 & 0.000 & 0.005 \\
		$\Gamma_{166}$ & 0.005 & 0.211 & 0.395 & 0.463 & 0.422 & 0.028 \\
		$\Gamma_{211}$ & 0.000 & 0.000 & 0.000 & 0.001 & 0.002 & 0.580 \\
		$\Gamma_{325}$ & 0.000 & 0.000 & 0.000 & 0.000 & 0.001 & 0.020 \\
		$\Gamma_{355}$ & 0.000 & 0.000 & 0.000 & 0.002 & 0.001 & 0.015 \\
		$\Gamma_{413}$ & 0.000 & 0.000 & 0.000 & 0.000 & 0.001 & 0.007 \\
		$\Gamma_{534}$ & 0.000 & 0.000 & 0.000 & 0.000 & 0.000 & 0.000 \\
		$\Gamma_{611}$ & 0.000 & 0.000 & 0.000 & 0.000 & 0.000 & 0.012 \\  
    \end{tabular}
    \end{small}
\begin{tiny}
\begin{align*}
&\Sigma_1 =
\begin{bmatrix}
      -0.000 &        0.000 &        0.000 &        0.000 &        0.000 &        0.000 \\ 
      -0.000 &       -0.000 &        0.000 &        0.000 &        0.000 &        0.010 \\ 
      -0.000 &        0.000 &       -0.000 &       -0.000 &        0.000 &       -0.000 \\ 
      -0.000 &        0.000 &        0.000 &        0.000 &        0.000 &       -0.000 \\ 
      -0.000 &        0.000 &       -0.000 &       -0.000 &        0.000 &       -0.010 \\ 
       0.000 &        0.000 &       -0.000 &       -0.000 &        0.000 &       -0.010 \\ 
\end{bmatrix} 
&\Sigma_2 =
\begin{bmatrix}
      -0.000 &        0.000 &        0.000 &        0.000 &        0.000 &       -0.000 \\ 
      -0.000 &       -0.001 &        0.000 &        0.001 &        0.001 &        0.733 \\ 
      -0.001 &        0.000 &       -0.001 &       -0.000 &        0.000 &        0.000 \\ 
      -0.000 &        0.000 &        0.000 &        0.000 &        0.000 &        0.000 \\ 
      -0.000 &        0.000 &       -0.000 &       -0.000 &        0.000 &       -0.423 \\ 
       0.000 &        0.000 &       -0.000 &       -0.000 &        0.000 &       -0.423 \\ 
\end{bmatrix} \\
&\Sigma_3 =
\begin{bmatrix}
      -0.003 &        0.000 &        0.000 &        0.004 &        0.000 &        0.000 \\ 
      -0.004 &       -0.018 &        0.000 &        0.017 &        0.093 &        5.283 \\ 
      -0.007 &        0.000 &       -0.010 &       -0.004 &        0.000 &       -0.000 \\ 
      -0.004 &        0.000 &        0.000 &        0.003 &        0.000 &       -0.000 \\ 
      -0.001 &        0.000 &       -0.000 &       -0.000 &        0.000 &       -0.837 \\ 
       0.000 &        0.000 &       -0.000 &       -0.000 &        0.000 &       -0.791 \\ 
\end{bmatrix}
&\Sigma_4 =
\begin{bmatrix}
      -0.028 &        0.000 &        0.000 &        0.043 &        0.000 &       -0.000 \\ 
      -0.328 &        0.247 &        0.012 &        1.226 &       11.951 &       46.268 \\ 
      -0.085 &        0.000 &       -0.111 &       -0.038 &        0.000 &        0.000 \\ 
      -0.041 &        0.000 &        0.000 &        0.035 &        0.000 &        0.000 \\ 
      -0.007 &        0.000 &       -0.000 &       -0.004 &        0.000 &       -0.975 \\ 
       0.000 &        0.000 &       -0.000 &       -0.000 &        0.000 &       -0.933 \\ 
\end{bmatrix} \\
&\Sigma_5 =
\begin{bmatrix}
      -0.028 &        0.000 &        0.000 &        0.043 &        0.000 &       -0.000 \\ 
      -0.018 &       -0.149 &        0.001 &        0.069 &        0.093 &        0.501 \\ 
      -0.085 &        0.000 &       -0.111 &       -0.038 &        0.000 &        0.000 \\ 
      -0.041 &        0.000 &        0.000 &        0.035 &        0.000 &        0.000 \\ 
      -0.012 &        0.000 &       -0.000 &       -0.007 &        0.000 &       -0.308 \\ 
       0.000 &        0.000 &       -0.000 &       -0.001 &        0.000 &       -0.859 \\ 
\end{bmatrix} 
&\Sigma_6 =
\begin{bmatrix}
       0.125 &        0.000 &        0.011 &        0.000 &        0.000 &       -0.000 \\ 
       2.104 &        0.495 &        0.059 &        0.002 &        0.873 &        0.001 \\ 
      -1.019 &        0.000 &        0.565 &       -0.005 &        0.000 &        0.000 \\ 
      -0.999 &        0.000 &        0.580 &        0.020 &        0.000 &        0.000 \\ 
      -2.206 &        0.000 &       -0.029 &       -0.001 &        0.000 &       -0.002 \\ 
      -0.303 &        0.000 &       -0.027 &       -0.001 &        0.000 &       -0.511 \\ 
\end{bmatrix}
\end{align*}
\end{tiny}
  \end{center}
\end{table}

In Case \#1, an initial guess very close to the solution is chosen, with a relative error of $-10^{-5}$ on the six nonlinear unknowns. NR's algorithm converges in just three iterations. As expected, the values of $\alpha_i$ and $\Gamma_{ijk}$ and $\sigma_{jj}$ are all below 0.01, indicating an excellent initial guess. 

In Case \#2, an initial guess with a $-0.1\%$ error is chosen for the six nonlinear unknowns. NR's algorithm converges in 5 iterations. The maximum $\alpha_i$ is $\alpha_1 = 0.224$, indicating that higher-order terms play some role in the first equation, while the maximum (and only significant non-zero) $\Gamma_{ijk}$ is $\Gamma_{166} = 0.211$, hence $\alpha = 0.224$ and $\beta = 0.211$. In this case, the sufficient condition of Theorem \ref{th:SufficientCondition} applies with $\alpha+\beta = 0.435$, which guarantees a reduction of the nonlinear residual after the first iteration with respect to the nonlinear residual computed with the initial guess of a factor about 2 or more, which is consistent with the relatively fast convergence of the algorithm.

The values $\alpha_1 = 0.224$ and $\Gamma_{166} = 0.211$ indicate that the only equation responsible for some non-negligible nonlinear behaviour in the NR iteration is the first one, and the responsible unknown variable is the sixth one, i.e. $p_i$. This is confirmed by matrix $\Sigma$, where the only significantly non-zero value of $\sigma_{jj}$ is the sixth one. 

In Case \#3, an initial guess with a $-1\%$ error is chosen for the six nonlinear unknowns. NR's algorithm fails after the first iteration, because the value of $p_i$ causes the argument of the square root in the first equation to become negative; hence, a value  $\lambda = 0.49 < 1$ must be taken to compute the $\alpha_i$ indicators. The largest values of the indicators, namely $\alpha_1 = 0.678$, $\Gamma_{166} = 0.395$, and $\abs{\sigma_{66}} = 0.791$ all indicate a problem with the initial guess of the sixth unknown $p_i$. Since the increment of $p_i$ in the first iteration is positive, one can try to improve the initial value of $p_i$ by increasing it, e.g. by halving the initial value of the term $p_s - p_i$ that appears under square root, which means $p_i = 2.1994$; this causes NR's algorithm to converge in 4 iterations.

In Case \#4, an initial guess with a $-10\%$ error is chosen for the six nonlinear unknowns. The situation is similar to the previous case, with somewhat higher values of the larger indicators $\alpha_1 = 1.316$, $\Gamma_{166} = 0.463$, and $\abs{\sigma_{66}} = 0.933$. As in the previous case, one could find a value that halves $p_s - p_i$, i.e. $p_i = 2.0905$, which however still causes convergence failure. Further halving $p_s - p_i$ leads to Case \#5, still with no convergence; the analysis of the largest $\alpha_i$, $\Gamma_{ijk}$ and $\abs{\sigma_{jj}}$ clearly indicates that the sixth unknown is still to blame, and should be further increased. By repetitively halving  $p_s - p_i$, the values of those indicators are progressively reduced, until convergence is achieved in 5 iterations with the initial guess $p_i$ = 2.1976. 

Note how the criteria introduced in the previous section clearly indicate how it is not necessary to change the initial guesses of the other five nonlinear unknowns in order to eventually reach convergence, despite the fact that their relative error is the same as the sixth one. This is due to the weaker nonlinearity of the equations involving them, which is reflected in the much lower values of the corresponding $\alpha_i$, $\Gamma_{ijk}$ and $\abs{\sigma_{jj}}$ factors. 

Finally, in Case \#6 the initial guess of all variables except the first is taken very close to the solution as in Case \#2; however, it is assumed that the initial guess of the first unknown $f$ is wrong by a factor 3, due to some gross mistake. NR's algorithm fails after the first iteration due to a negative square root argument. The largest gamma value is $\Gamma_{211} = 0.58$, indicating a problem with the first unknown in the second equation. The three largest sigma values are $\abs{\sigma_{33}} = 0.565$,  $\abs{\sigma_{66}} = 0.511$, and $\abs{\sigma_{33}} = 0.495$, which are very close to each other, potentially indicating a problem with the second, third, and sixth unknown. In this case, the gamma values turn out to be more selective, as they clearly point to the wrong initial guess of $f$, which can be then spotted immediately and corrected. Having done that, the algorithm converges in a few iterations.

\subsection{DC circuit example}
Consider an electrical DC circuit, where the series connection of $N$ resistors and one diode is connected to an ideal voltage source, which provides a certain power $P$. The system is described by the following set of implicit equations:
\begin{align}
&i - \left( i_s e^{v_d/v_t} - 1 \right) = 0\\
&vi - P = 0\\
&v - \sum_{j=1}^N v_j - v_d = 0\\
&v_j - R i = 0
\end{align}
where $i_s, v_t, P, R$ are known parameters, $x = \begin{bmatrix}i & v_d & v & v_1 & v_2 & \cdots & v_N\end{bmatrix}'$,  $w = \begin{bmatrix}i & v_d & v \end{bmatrix}'$, $z = \begin{bmatrix} v_1 & v_2 & \cdots & v_N\end{bmatrix}'$.


%

Taking $i_s = 6.9144\cdot 10^{-13}, v_t = 25\cdot 10^{-3}, P = 10.7, R = 1, N = 10$, the system has an exact solution $\bar{i} = 1, \bar{v}_d = 0.7, \bar{v} = 10.7,  \bar{v}_j = 1$.

There are only three non-zero curvature factors for this problem, namely $\Gamma_{122}$, corresponding to $v_d$ in the first equation, and $\Gamma_{213} = \Gamma_{231}$, corresponding to $i$ and $v$ in the second equation.

According to Criterion \ref{cr:InitialGuesses}, in order to ensure fast and
reliable convergence, accurate initial guesses should be provided for the
unknowns $w$, while $z_0$ can be safely taken to be zero, since other choice
leads to exactly the same results after each iteration. In this specific case,
this means that only 3 variables out of 13 require to be properly initialized.
Some experiments indeed confirmed that the results after the first iteration are not
affected at all from the values of $z_0$.

Table \ref{tab:CircuitConvergence} reports the initial guesses $w_0$, the number of NR iterations, and all the relevant $\alpha_i$, $\Gamma_{ijk}$, and $\Sigma$ indicators, corresponding to different choices of $w_0$.

\begin{table}[hbt]
  \begin{center}
    \caption{Convergence analysis of the DC circuit test case}
    \label{tab:CircuitConvergence}
    \begin{small}
    \begin{tabular}{c|c|c|c|c|c}    
       Var           & \#1        & \#2      & \#3               & \#4                 & \#5    \\
       \hline
       $i$           & 0.99999    &  0.99    & 0.9               & 0.8                 & 0.25   \\
       $v_d$         & 0.699993   &  0.693   & 0.63              & 0.56                & 0.693  \\
       $v$           & 10.699893  & 10.593   & 9.63              & 8.56                & 2.675  \\
       $N_{iter}$    & 2          &  4       & 18		         & --                  & 7      \\
       $\alpha_{1}$   & 0.000 & 0.020 & $1.31\cdot 10^5$ & $1.18\cdot 10^{88}$ & 0.071 \\
       $\Gamma_{122}$ & 0.000 & 0.168 &    3.497 &   21.116 & 0.067 \\
	   $\Gamma_{213}$ & 0.000 & 0.002 & 0.029 & 0.014 & 0.958 \\
    \end{tabular}
    \end{small}
\begin{align*}
&\Sigma_1 =
\begin{bmatrix}
      -0.000 &        0.000 &       -0.000 \\
      -0.000 &       -0.000 &       -0.000 \\
      -0.000 &       -0.000 &       -0.000 \\
\end{bmatrix} 
&\Sigma_2 =
\begin{bmatrix}
      -0.005 &        0.013 &       -0.005 \\
      -0.000 &       -0.323 &       -0.000 \\
      -0.005 &       -0.012 &       -0.005 \\
\end{bmatrix} \\
&\Sigma_3 =
\begin{bmatrix}
      -0.068 &        3.054 &       -0.068 \\
      -0.007 &      -14.993 &       -0.007 \\
      -0.050 &       -2.299 &       -0.050 \\
\end{bmatrix} 
&\Sigma_4 =
\begin{bmatrix}
       0.229 &     1934.463 &        0.229 \\
       0.006 &     -158.105 &        0.006 \\
       0.016 &      -85.085 &        0.016 \\
\end{bmatrix} \\
&\Sigma_5 =
\begin{bmatrix}
      -3.796 &        0.002 &       -3.796 \\
      -5.162 &       -1.856 &       -5.162 \\
      -3.699 &       -0.002 &       -3.699 \\
\end{bmatrix}
\end{align*}
    
  \end{center}
\end{table}

Note that $\alpha_2 = 0$ in all cases; this is due to the fact that the residual of the second equation is a second-order polynomial, so obviously its Taylor expansion lacks terms above second order.

In Case \#1, an initial guess very close to the solution is chosen, with a relative error of $-10^{-5}$ on the three nonlinear unknowns. NR's algorithm converges in just two iterations. As expected, the values of $\alpha_i$, $\Gamma_{ijk}$ and $\abs{\sigma_{jj}}$ are all below 0.001, indicating an excellent initial guess. 

In Case \#2, an initial guess with a $-1\%$ error is chosen for the three nonlinear unknowns. NR's algorithm converges in 4 iterations. The maximum $\alpha_i$ is $\alpha_1 = 0.02$, indicating that higher-order terms play a negligible role in the first equation, while the maximum $\Gamma_{ijk}$ is $\Gamma_{122} = 0.168$. Here, $\alpha = 0.02$ and $\beta = 0.170$. In this case, the sufficient condition of Theorem \ref{th:SufficientCondition} applies with $\alpha+\beta = 0.190$, which guarantees a reduction of the nonlinear residual after the first iteration with respect to the nonlinear residual computed with the initial guess of a factor about 5, which is consistent with the very fast convergence of the algorithm. The values $\Gamma_{122} = 0.168$ and $\abs{\sigma_{22}} = 0.323$ indicate that the second nonlinear unknown $v_d$ is solely responsible for some non-negligible nonlinear behaviour in the first equation, which is however small enough not to cause problems. This is obviously due to the strongly exponential behaviour of the diode equation with respect to $v_d$.

In Case \#3, an initial guess with a $-10\%$ error is chosen for the three
nonlinear unknowns. NR's algorithm requires 18 iterations to converge. The
large value $\alpha_1 = 1.31\cdot 10^5$ clearly indicates that the cause of this
difficult convergence is an error in the initial guess of only nonlinear
variable appearing in the first equation, i.e., $v_d$. The other two larger indicators $\Gamma_{122} = 3.497$ and $\abs{\sigma_{22}} = 14.993$ confirm this finding, as they both point to the second unknown $v_d$. The positive increment of $v_d$ after one iteration suggests to increase its initial guess; in fact, increasing it by 0.05 to $v_d = 0.68$ reduces the number of iterations to 8, while increasing it by 0.10 to $v_d$ = 0.73 reduces the number of iterations to 6.
Note that it is not necessary to worry about the other initial guesses to
substantially improve the convergence performance, as suggested by the much lower values of all the other indicators.

In Case \#4, an initial guess with a $-20\%$ error is chosen for the three nonlinear unknowns. NR's algorithm fails at the second iteration because of badly conditioned Jacobian. The value $\alpha_1 = 1.18 \cdot 10^{88}$ reveals a huge contribution of the higher-order terms in the first equation - this is caused by the exponential in the first equation, which has a scaling factor $v_t = 0.025$, meaning that errors in $v_d$ significantly larger than the scaling factor have a dramatic impact on the behaviour of the equation residual. The result is a large overshoot of $w_{1,2} - w_{0,2}$ that brings $v_d$ so far away from the solution that the following iteration is not even possible to compute. This is also confirmed by the largest values $\Gamma_{122} = 21.116$ and $\abs{\sigma_{22}} = 158.105$.

As in the previous case, the positive sign of the increment of $v_d$ suggests to increase its initial guess, while the low values of $\alpha_2$, $\Gamma_{213}$, $\sigma_{11}$, and $\sigma_{33}$ indicate that the initial guesses of $v$ and $i$ are close enough to the solution so as to not be a problem. An increase by 0.05 to $v_d = 0.61$ leads to convergence in 37 steps; a further increase by 0.05 to $v_d = 0.66$ leads to convergence in 8 steps. Throughout these steps, all the indicators pertaining to the other two nonlinear variables remain well below unity.

In Case \#5, the initial guess of $v_d$ has a small error of $-1\%$, while the initial guess of the other two variables has a much larger error of $-75\%$. NR's algorithm converges in 7 steps, which is acceptable, but could be probably improved. Also in this case, Criterion \ref{cr:Ranking} provides the correct diagnosis, as the largest indicator values $\Gamma_{213} = 0.958$, $\abs{\sigma_{11}} = 3.796$, and $\abs{\sigma_{33}} = 3.699$ clearly indicate that there is something wrong with the first and third unknown in the second equation.

\subsection{AC distribution system}
\label{sec:ACSystem}
The last example case is the power flow of an AC balanced three-phase grid with $N \times N$ nodes, $N$ being an even integer. Three-phase voltages and currents are described by complex phasors. The even nodes are connected to power generators, that prescribe the voltage magnitude $V_g$ and the amount of injected  active electric power $P$. The odd nodes are connected to resistive loads, which absorb the same active power $P$ and no reactive power, and have a voltage modulus $V_l$. The nodes are connected by a square grid of purely inductive transmission lines with impedance $Z = jX$ and an admittance $Y = 1/jX$, where $j = \sqrt{-1}$.

The system is balanced, because there are $N^2/2$ generators collectively injecting an active power $PN^2/2$, and an equal number of loads collectively absorbing the same amount of active power, while no active power is lost in the transmission lines which have zero resistance. In order to obtain a well-posed power flow problem, the generator at node $(g,g)$, $g = N/2$ is substituted by a slack node, that fixes the node voltage phasor to the known complex value $V_g + j\cdot 0$.

The system is described by the following equations:
\begin{align}
& \abs{v_{i,k}} = V_g & i,k = 1, \cdots, N \quad i,k \neq g \quad i+k \; \mathrm{even}\\
& \mathrm{Re}(v_{i,k} \bar{i}^n_{i,k}) = - P  & i,k = 1, \cdots, N \quad i,k \neq g \quad i+k \; \mathrm{even}\\
& v_{i,k} \bar{i}^n_{i,k} = P  & i,k = 1, \cdots, N \quad i+k \; \mathrm{odd}\\
& v_{g,g} = V_g\\
& i^h_{i,k} = Y(v_{i,k} - v_{i+1,k}) & i = 1, \cdots, N-1 \quad k = 1, \cdots, N \\
& i^v_{i,k} = Y(v_{i,k} - v_{i,k+1}) & i = 1, \cdots, N \quad k = 1, \cdots, N-1 \\
& i^n_{1,1} + i^h_{1,1} + i^v_{1,1} = 0 \\
& i^n_{N,1} - i^h_{N-1,1} + i^v_{N,1} = 0 \\
& i^n_{1,N} - i^v_{1,N-1} + i^h_{1,N} = 0 \\
& i^n_{N,N} - i^h_{N-1,N-1} - i^v_{N,N-1} = 0 \\
& i^n_{i,1} + i^h_{i,1} + i^v_{i,1} - i^h_{i-1,1} = 0 & i = 2, \cdots, N-1\\
& i^n_{1,k} + i^h_{1,k} + i^v_{1,k} - i^v_{1,k-1} = 0 & k = 2, \cdots, N-1\\
& i^n_{N,k} + i^v_{N,k} - i^h_{N-1,k} - i^v_{N,k-1} = 0, & k = 2, \cdots, N-1\\
& i^n_{i,N} + i^h_{i,N} - i^h_{i-1,N} - i^v_{i,N-1} = 0, & i = 2, \cdots, N-1\\
& i^n_{i,k} + i^h_{i,k} + i^v_{i,k} -i^h_{i-1,k} - i^v_{i,k-1} = 0, & i,k = 2, \cdots, N-1,
\end{align}
where
\begin{itemize}
\item $v_{i,k}$ are voltages of nodes $(i,k)$, $i,k = 1, \cdots, N$;
\item $i^n_{i,k}$ are the currents leaving the $(i,k)$ node and entering either the generator or load, $i,k = 1, \cdots, N$;
\item $i^h_{i,k}$ are the currents flowing through the horizontal lines connecting nodes $(i,k)$ and $(i+1,k)$, $i = 1, \cdots, N-1$, $k = 1, \cdots, N$;
\item $i^v_{i,k}$ are the currents flowing through the vertical lines connecting nodes $(i,k)$ and $(i,k+1)$, $i = 1, \cdots, N$, $k = 1, \cdots, N-1$,
\end{itemize}
for a total of $4N^2 - 2N$ complex unknowns. The complex equations can be split into their real and imaginary parts, involving  $8N^2 - 4N$ real unknowns, which are the real and imaginary coefficients of the complex unknowns. The bar denotes the complex conjugate operator. All variables and parameters are given in per unit, using the active power consumed by the load as base power, and the voltage modulus of the load as base voltage, so that the system is well scaled by construction.

More specifically, the first two equations correspond to the power generators; the following equation corresponds to the loads; the following equation corresponds to the slack node; the following two equations are Ohm's law applied to the horizontal and vertical lines; the following four equations are Kirchhoff's current laws applied to the nodes in the top left, top right, bottom left and bottom right corners of the grid; the following four equations are Kirchhoff's current laws applied to the nodes of the top, left, right, and bottom edges of the grid; the last equation is Kirchhoff's law applied to the internal nodes of the grid. Overall, there are $2N^2 - 2$ nonlinear equations and  $6N^2 - 4N + 2$ linear equations in the system.

In the limit case of an infinite grid, all generator node voltages must be equal and all load voltages must be equal, because of symmetry. Each generator feeds four adjacent loads by means of four lines in parallel, and each load is fed by four adjacent generators, by means of four equal lines in parallel. Assuming the generator voltages are real (i.e., they have zero phase) and that those voltages are set to get a unit voltage magnitude of the loads, it is thus straightforward to compute the exact solution in the infinite grid case:
\begin{align}
& P = 1 \\
& V_g = \sqrt{1+\frac{X^2}{16}} \\
& V_l = \frac{V_g}{1 + j \frac{X}{4}}\\
& v_{i,k} = V_g & i+k \; \mathrm{even} \\
& v_{i,k} = V_l & i+k \; \mathrm{odd}\\
& i^n_{i,k} = -V_l &  i+k \; \mathrm{even} \\
& i^n_{i,k} = V_l & i+k \; \mathrm{odd} \\
& i^h_{i,k} = Y(v_{i,k} - v_{i+1,k}) \\
& i^v_{i,k} = Y(v_{i,k} - v_{i,k+1})
\end{align}

In the case of a finite $N \times N$ grid, the solution is similar to that of the infinite grid in the interior nodes which are far from the edges, while border effects warp the solution when getting close to the edges. In particular, the two generators at the grid corners have to feed their active power $P$ through two lines only, instead of the four lines available at the center of the grid, while the generators along the edges have to feed their active power through three lines. As a consequence, their voltage must be significantly higher than in the ideal infinite grid solution, because they must overcome a higher impedance (two or three lines in parallel instead of four) to feed the neighbouring loads.

The results of a several test experiments on a $20 \times 20$ grid, with 798 nonlinear and 2322 linear unknowns, are now presented. Note that the original problem involves complex numbers and equations, which were split into real and imaginary part for convenience in writing the code to solve the problem and compute the indicators. In case an indicator points to the real or imaginary part of the initial guess of complex unknown, the entire complex value should be considered. More detailed results of the analysis are reported in the supplementary material.

In Test \#1, the infinite grid solution is used as an initial guess for the problem; NR's method converges in 7 iterations. Given the analysis reported above, it is expected that the criteria presented in the previous section show that the initial guesses for the variables at the corners and edges of the grid are the farthest from the solution. In fact, the two largest values of $\Gamma_{ijk}$, around 0.6, the two largest values of $\alpha_i$, around 0.1, and the two largest values of $\abs{\sigma_{jj}}$, around 3, all point to the voltages of the two generators at the grid corners. Subsequent indicators in the ranking point to the voltages of the two loads at the grid corner, and to some generator voltages along the edges of the grid, close to the corners.

In subsequent tests, the initial guess is obtained by applying changes to some selected variables, on top of the infinite grid solution, so that they get farther from the solution. The objective of the analysis is to ensure that Criterion \ref{cr:Ranking} always succeeds at pinpointing those variables that were changed, possibly causing NR's algorithm to fail, and that fixing them by ranking order can eventually restore the convergence.

In Test \#2, the initial guess of $i^n_{5,1}$ was set to 1/10 of the exact solution. Apparently, such a large change does not have a major impact on convergence, which is achieved in 6 iterations, one less than in the previous case. In fact, the highest ranking indicators still point to voltage variables corresponding to nodes on the corners and edges of the grid.

In Test \#3, the initial guess of $i^n_{5,1}$ was set to 1/10 of the exact solution, and the initial guess of the corresponding voltage $v_{5,1}$ was set to 1/2 of the exact solution. NR's method now converges after 7 iterations. $v_{5,1}$ now shows up as first in the ranking of $\alpha_i$ with a value of 0.302, and second in the ranking of $\Gamma_{ijk}$ with a value of 0.588. Subsequent values in the ranking still point to voltages of the corners and edges of the grid.

If the initial guess of $v_{5,1}$ is reduced to 1/10 of the exact solution, as in Test \#4, NR's method fails to converge. Variable $v_{5,1}$ ranks first with $\alpha_i = 4690$ and with $\Gamma_{ijk} = 4370$, while both $v_{5,1}$ and $i^n_{5,1}$ rank first with similar values around $\abs{\sigma_{jj}} = 10000$. One then need to get the initial values of $v_{5,1}$ and, possibly, of $i^n_{5,1}$, closer to the solution in order to achieve convergence.

Reversing the situation of Test \#3 brings to Test \#5, with $v_{5,1}$ initialized at 1/10 of the exact solution and $i^n_{5,1}$ at 1/2 of the exact solution. In this case, NR's method does not converge and the culprit is very clearly indicated in $v_{5,1}$ by the values $\alpha_i = 123$, $\Gamma_{ijk}$ = 101, and $\abs{\sigma_{jj}} = 1070$, which are by far the top ranking ones.

Test \#6 is similar to Test \#2, with the initial guess of $v_{5,1}$ set to 1/100 of the exact value; convergence is achieved in 7 iterations. Variable $v_{5,1}$ is clearly indicated as responsible of the less-than-ideal convergence by the top ranking indicators (each in its category) $\alpha_i = 0.996$, $\Gamma_{ijk} = 1.43$, and $\sigma_{jj} = 37$.

In the subsequent Test \#7, several initial guess were altered: $v_{5,1}$ to 1/2 of the exact solution; $v_{13,17}$, $v_{13,18}$, $v_{14,18}$, $v_{12,20}$, $v_{13,20}$, $i^n_{5,1}$, $i^n_{13,17}$, $i^n_{13,18}$, $i^n_{14,18}$, and $i^n_{13,20}$ to 1/10 of the exact solution. NR's method does not converge. Variable $v_{12,20}$ ranks first with $\abs{\sigma_{jj}} = 174$ and with a factor 4 margin on the next smaller $\sigma_{jj} = 46.3$ value. It also ranks first with $\alpha_i = 3.49$, although with a tiny margin over the next $\alpha_i = 3.4$ value, and ranks sixth with $\Gamma_{ijk} = 3.8$, the largest one having $\Gamma_{ijk} = 5.41$, which is not much larger. One can then first try to improve $v_{12,20}$.

Setting the initial guess of $v_{12,20}$ to 80\% of the exact value brings to Test \#8, where $v_{13,20}$ ranks first with $\Gamma_{ijk}$ = 3.12 and with $\abs{\sigma_{jj}} = 72$. If $v_{13,20}$ is then also brought to 80\% of the exact solution, Test \#9 is obtained, which still does not converge. The highest-ranking $\alpha_i = 4.25$ and $\Gamma_{ijk} = 4.62$ point to $v_{14,18}$.

If one now changes $v_{14,18}$ to 80\% of the exact value, Test \#10 is obtained, which still doesn't converge and shows the highest ranking $\alpha_i = 1.6$ and $\Gamma_{ijk} = 1.79$ for $v_{13,17}$. Further setting the initial guess for this variable to 80 \% of the exact solution brings us to Test \#11, which does then converge in 13 iterations.

This battery of tests clearly demonstrates that the proposed indicators succeed in pointing out the most critical initial guesses, which are known a priori in this experimental setting, and that the strategy of fixing initial guesses one by one according to the ranking of the indicators is effective to eventually achieve convergence.

\section{Conclusion}
\label{sec:Conclusion}
In this paper, new theorems were presented concerning the choice of initial guess values for NR's method, when solving generic systems of nonlinear or mixed linear and nonlinear equations. Based on these theorems, four criteria were proposed to help choosing or improving initial guess values for NR's method, in order to achieve convergence.

Criterion \ref{cr:InitialGuesses}, which is rigorous and only based on structural properties of the system of equations, suggests to only care about the initial guess of the subset of variables that influence the Jacobian of the nonlinear system; all other variables can be initially set to zero without any consequence on the convergence of NR's algorithm. 

Criterion \ref{cr:Ranking} can be used in case of convergence failure of the NR solver, to identify those variables whose initial guesses are most likely the cause of the failure, also suggesting how to improve them to eventually achieve successful convergence of the iterative solution process.

Criterion \ref{cr:Ranking} makes use of the indicators $\alpha_i$, $\Gamma_{ijk}$ and $\abs{\sigma_{jj}}$, defined in Sections \ref{sec:Method} and \ref{sec:Discussion}, which provide information about the second- and higher-order behaviour of the function residuals around the chosen initial guess. The computation of these indicators require the Jacobian and Hessian matrices of the equation residuals; these can be computed analytically by symbolic differentiation, which is the standard approach taken by EOOLTs, or by numerical differentiation. 

The $\alpha_i$ and $\abs{\sigma_{jj}}$ indicators are invariant with respect to the scaling of the problem, i.e., they are not affected by the choice of measurement units of the variables of the involved variables, which is an important property for equations coming from physical system modelling. The $\Gamma_{ijk}$ indicators instead require proper scaling of the equation residuals, in order to obtain meaningful results.

In the context of EOOLTs, Criterion \ref{cr:InitialGuesses} can be used by
developers of reusable model libraries, to identify which variables actually
need good initial guesses, and thus to provide the proper infrastructure to do
so, e.g. by providing ad-hoc parameters to be set by the end users. Criterion \ref{cr:Ranking} can instead be used by simulation tool developers to
provide meaningful diagnostic information in case of NR solver failures,
guiding the end user towards the successful solution of the problem by means of a suitable graphical user interface.

Criterion \ref{cr:Ranking} was then successfully demonstrated in three exemplary cases discussed in Section \ref{sec:Examples}: a simple thermo-hydraulic system, an electrical DC circuit involving a diode and several resistors, and a large electrical AC power flow problem. Several tests were set up for each example case, in which initial guesses were suitably changed with respect to the exact solution. In all cases, the proposed criterion successfully identified the few initial guesses that need to be improved, eventually leading to convergence of NR's algorithm.

Even though this research was originally motivated by the need of good diagnostic tools for the initialization of equation-based, object-oriented dynamic models, the results presented in this paper are of course not limited to that case, but rather have a very broad applicability to any kind of problem that requires the solution of implicit nonlinear (or mixed linear and nonlinear) systems of equations by means of NR's algorithm.

\section{Acknowledgements}
The authors are grateful to colleagues Luca Bonaventura and Gianni Ferretti of Politecnico di Milano for their constructive comments on the first draft of the paper.

This research did not receive any specific grant from funding agencies in the public, commercial, or not-for-profit sectors.

\bibliographystyle{abbrv}
\bibliography{oom,PublicationsCasella,paper}

\end{document}


\maketitle
\newpage
Example 1: Using solution of the infinite grid as default guess values.

\vspace{-0.5cm}
\begin{table}[H]
	\begin{center}
	\caption{Sorted $\Gamma$ values with corresponding equation and variables.}
	\begin{small}
	\begin{tabular}{c||l|l|l}
		$\Gamma_{ijk}$ & equation $i$ & variable $j$ & variable $k$\\
	\hline
       0.669 & $\abs{v_{1,1}} = V_g$ & $\mathrm{Im}(v_{1,1}) =     0$ & $\mathrm{Im}(v_{1,1}) =     0$ \\
       0.626 & $\abs{v_{20,20}} = V_g$ & $\mathrm{Im}(v_{20,20}) =     0$ & $\mathrm{Im}(v_{20,20}) =     0$ \\
       0.374 & $\mathrm{Re}(v_{20,1} \bar{i}^n_{20,1}) = \mathrm{Re}(P)$ & $\mathrm{Im}(v_{20,1}) = -0.243$ & $\mathrm{Im}(i^n_{20,1}) = -0.243$ \\
       0.374 & $\mathrm{Re}(v_{1,20} \bar{i}^n_{1,20}) = \mathrm{Re}(P)$ & $\mathrm{Im}(v_{1,20}) = -0.243$ & $\mathrm{Im}(i^n_{1,20}) = -0.243$ \\
       0.365 & $\mathrm{Re}(v_{1,1} \bar{i}^n_{1,1})) = \mathrm{Re}(-P)$ & $\mathrm{Im}(v_{1,1}) =     0$ & $\mathrm{Im}(i^n_{1,1}) = 0.243$ \\
        0.34 & $\mathrm{Re}(v_{20,20} \bar{i}^n_{20,20})) = \mathrm{Re}(-P)$ & $\mathrm{Im}(v_{20,20}) =     0$ & $\mathrm{Im}(i^n_{20,20}) = 0.243$ \\
	\end{tabular}
	\end{small}
	\end{center}
\end{table}

\vspace{-1cm}
\begin{table}[H]
	\begin{center}
	\caption{Sorted $\alpha$ and $\Sigma$ values with corresponding equations and variables, respectively.}
	\begin{small}
	\begin{tabular}{c c c}
	\begin{tabular}{c||l}
		$\alpha_i$ & equation $i$ \\
	\hline
       0.128 & $\abs{v_{1,1}} = V_g$ \\
       0.114 & $\abs{v_{20,20}} = V_g$ \\
       0.035 & $\abs{v_{1,3}} = V_g$ \\
       0.035 & $\abs{v_{3,1}} = V_g$ \\
      0.0299 & $\abs{v_{18,20}} = V_g$ \\
      0.0299 & $\abs{v_{20,18}} = V_g$ \\
	\end{tabular} & &
	\begin{tabular}{c||l}
		 $\sigma_{j,j}$ & variable $j$ \\
	\hline
       -3.06 & $\mathrm{Im}(v_{20,1}) = -0.243$ \\
       -3.06 & $\mathrm{Im}(v_{1,20}) = -0.243$ \\
       -2.03 & $\mathrm{Im}(i^n_{1,20}) = -0.243$ \\
       -2.03 & $\mathrm{Im}(i^n_{20,1}) = -0.243$ \\
        1.93 & $\mathrm{Im}(v_{20,2}) =     0$ \\
        1.93 & $\mathrm{Im}(v_{2,20}) =     0$ \\
	\end{tabular}
	\end{tabular}
	\end{small}
	\end{center}
\end{table}

\vspace{-.5cm}
SUCCESS: Newton-Raphson algorithm converged in 7 iterations!

\newpage
Example 2: Changing some guess values (see table below).

\vspace{-0.5cm}
\begin{table}[H]
	\begin{center}
	\caption{Critical initial guess values compared to exact solution.}
	\begin{small}
	\begin{tabular}{l||c|c}
		variable & initial guess & exact solution \\
	\hline
		$i^n_{ 5, 1}$  & $ -0.1013 + i\;   0.0128$ & $ -1.0133 + i\;   0.1276$ \\
	\end{tabular}
	\end{small}
	\end{center}
\end{table}

\vspace{-1cm}
\begin{table}[H]
	\begin{center}
	\caption{Sorted $\Gamma$ values with corresponding equation and variables.}
	\begin{small}
	\begin{tabular}{c||l|l|l}
		$\Gamma_{ijk}$ & equation $i$ & variable $j$ & variable $k$\\
	\hline
        0.81 & $\abs{v_{20,20}} = V_g$ & $\mathrm{Im}(v_{20,20}) =     0$ & $\mathrm{Im}(v_{20,20}) =     0$ \\
        0.56 & $\abs{v_{1,1}} = V_g$ & $\mathrm{Im}(v_{1,1}) =     0$ & $\mathrm{Im}(v_{1,1}) =     0$ \\
       0.448 & $\mathrm{Re}(v_{20,20} \bar{i}^n_{20,20})) = \mathrm{Re}(-P)$ & $\mathrm{Im}(v_{20,20}) =     0$ & $\mathrm{Im}(i^n_{20,20}) = 0.243$ \\
       0.407 & $\mathrm{Re}(v_{19,19} \bar{i}^n_{19,19})) = \mathrm{Re}(-P)$ & $\mathrm{Im}(v_{19,19}) =     0$ & $\mathrm{Im}(i^n_{19,19}) = 0.243$ \\
       0.406 & $\abs{v_{18,20}} = V_g$ & $\mathrm{Im}(v_{18,20}) =     0$ & $\mathrm{Im}(v_{18,20}) =     0$ \\
       0.406 & $\abs{v_{20,18}} = V_g$ & $\mathrm{Im}(v_{20,18}) =     0$ & $\mathrm{Im}(v_{20,18}) =     0$ \\
	\end{tabular}
	\end{small}
	\end{center}
\end{table}

\vspace{-1cm}
\begin{table}[H]
	\begin{center}
	\caption{Sorted $\alpha$ and $\Sigma$ values with corresponding equations and variables, respectively.}
	\begin{small}
	\begin{tabular}{c c c}
	\begin{tabular}{c||l}
		$\alpha_i$ & equation $i$ \\
	\hline
       0.175 & $\abs{v_{20,20}} = V_g$ \\
      0.0943 & $\abs{v_{1,1}} = V_g$ \\
      0.0539 & $\abs{v_{18,20}} = V_g$ \\
      0.0538 & $\abs{v_{20,18}} = V_g$ \\
      0.0445 & $\abs{v_{19,19}} = V_g$ \\
       0.024 & $\abs{v_{1,3}} = V_g$ \\
	\end{tabular} & &
	\begin{tabular}{c||l}
		 $\sigma_{j,j}$ & variable $j$ \\
	\hline
       -2.65 & $\mathrm{Im}(v_{20,1}) = -0.243$ \\
       -2.45 & $\mathrm{Im}(v_{1,20}) = -0.243$ \\
       -1.72 & $\mathrm{Im}(i^n_{20,1}) = -0.243$ \\
        1.65 & $\mathrm{Im}(v_{20,19}) = -0.243$ \\
        1.65 & $\mathrm{Im}(v_{19,20}) = -0.243$ \\
        1.62 & $\mathrm{Im}(v_{19,1}) =     0$ \\
	\end{tabular}
	\end{tabular}
	\end{small}
	\end{center}
\end{table}

\vspace{-.5cm}
SUCCESS: Newton-Raphson algorithm converged in 6 iterations!

\newpage
Example 3: Changing some guess values (see table below).

\vspace{-0.5cm}
\begin{table}[H]
	\begin{center}
	\caption{Critical initial guess values compared to exact solution.}
	\begin{small}
	\begin{tabular}{l||c|c}
		variable & initial guess & exact solution \\
	\hline
		$v_{ 5, 1}$ & $  0.5058 + i\;   0.0987$ & $  1.0117 + i\;   0.1974$ \\
		$i^n_{ 5, 1}$  & $ -0.0101 + i\;   0.0013$ & $ -1.0133 + i\;   0.1276$ \\
	\end{tabular}
	\end{small}
	\end{center}
\end{table}

\vspace{-1cm}
\begin{table}[H]
	\begin{center}
	\caption{Sorted $\Gamma$ values with corresponding equation and variables.}
	\begin{small}
	\begin{tabular}{c||l|l|l}
		$\Gamma_{ijk}$ & equation $i$ & variable $j$ & variable $k$\\
	\hline
       0.593 & $\abs{v_{1,1}} = V_g$ & $\mathrm{Im}(v_{1,1}) =     0$ & $\mathrm{Im}(v_{1,1}) =     0$ \\
       0.588 & $\abs{v_{5,1}} = V_g$ & $\mathrm{Im}(v_{5,1}) = 0.0987$ & $\mathrm{Im}(v_{5,1}) = 0.0987$ \\
       0.445 & $\mathrm{Re}(v_{1,20} \bar{i}^n_{1,20}) = \mathrm{Re}(P)$ & $\mathrm{Im}(v_{1,20}) = -0.243$ & $\mathrm{Im}(i^n_{1,20}) = -0.243$ \\
       0.363 & $\mathrm{Re}(v_{20,1} \bar{i}^n_{20,1}) = \mathrm{Re}(P)$ & $\mathrm{Im}(v_{20,1}) = -0.243$ & $\mathrm{Im}(i^n_{20,1}) = -0.243$ \\
       0.352 & $\abs{v_{3,1}} = V_g$ & $\mathrm{Im}(v_{3,1}) =     0$ & $\mathrm{Im}(v_{3,1}) =     0$ \\
       0.336 & $\mathrm{Re}(v_{1,1} \bar{i}^n_{1,1})) = \mathrm{Re}(-P)$ & $\mathrm{Im}(v_{1,1}) =     0$ & $\mathrm{Im}(i^n_{1,1}) = 0.243$ \\
	\end{tabular}
	\end{small}
	\end{center}
\end{table}

\vspace{-1cm}
\begin{table}[H]
	\begin{center}
	\caption{Sorted $\alpha$ and $\Sigma$ values with corresponding equations and variables, respectively.}
	\begin{small}
	\begin{tabular}{c c c}
	\begin{tabular}{c||l}
		$\alpha_i$ & equation $i$ \\
	\hline
       0.302 & $\abs{v_{5,1}} = V_g$ \\
       0.161 & $\abs{v_{1,1}} = V_g$ \\
      0.0693 & $\abs{v_{3,1}} = V_g$ \\
      0.0534 & $\abs{v_{1,3}} = V_g$ \\
      0.0517 & $\abs{v_{2,2}} = V_g$ \\
      0.0307 & $\abs{v_{4,2}} = V_g$ \\
	\end{tabular} & &
	\begin{tabular}{c||l}
		 $\sigma_{j,j}$ & variable $j$ \\
	\hline
       -4.04 & $\mathrm{Im}(v_{1,20}) = -0.243$ \\
       -3.93 & $\mathrm{Im}(v_{20,1}) = -0.243$ \\
       -2.89 & $\mathrm{Im}(v_{5,1}) = 0.0987$ \\
        -2.7 & $\mathrm{Im}(i^n_{1,20}) = -0.243$ \\
       -2.65 & $\mathrm{Im}(i^n_{20,1}) = -0.243$ \\
        2.63 & $\mathrm{Im}(v_{2,20}) =     0$ \\
	\end{tabular}
	\end{tabular}
	\end{small}
	\end{center}
\end{table}

\vspace{-.5cm}
SUCCESS: Newton-Raphson algorithm converged in 7 iterations!

\newpage
Example 4: Changing some guess values (see table below).

\vspace{-0.5cm}
\begin{table}[H]
	\begin{center}
	\caption{Critical initial guess values compared to exact solution.}
	\begin{small}
	\begin{tabular}{l||c|c}
		variable & initial guess & exact solution \\
	\hline
		$v_{ 5, 1}$ & $  0.0101 + i\;   0.0020$ & $  1.0117 + i\;   0.1974$ \\
		$i^n_{ 5, 1}$  & $ -0.0101 + i\;   0.0013$ & $ -1.0133 + i\;   0.1276$ \\
	\end{tabular}
	\end{small}
	\end{center}
\end{table}

\vspace{-1cm}
\begin{table}[H]
	\begin{center}
	\caption{Sorted $\Gamma$ values with corresponding equation and variables.}
	\begin{small}
	\begin{tabular}{c||l|l|l}
		$\Gamma_{ijk}$ & equation $i$ & variable $j$ & variable $k$\\
	\hline
    4.37e+03 & $\abs{v_{5,1}} = V_g$ & $\mathrm{Im}(v_{5,1}) = 0.00197$ & $\mathrm{Im}(v_{5,1}) = 0.00197$ \\
         158 & $\abs{v_{5,1}} = V_g$ & $\mathrm{Re}(v_{5,1}) = 0.0101$ & $\mathrm{Im}(v_{5,1}) = 0.00197$ \\
          46 & $\mathrm{Re}(v_{5,1} \bar{i}^n_{5,1})) = \mathrm{Re}(-P)$ & $\mathrm{Im}(v_{5,1}) = 0.00197$ & $\mathrm{Im}(i^n_{5,1}) = 0.00128$ \\
        43.2 & $\abs{v_{20,20}} = V_g$ & $\mathrm{Im}(v_{20,20}) =     0$ & $\mathrm{Im}(v_{20,20}) =     0$ \\
        40.7 & $\mathrm{Re}(v_{19,19} \bar{i}^n_{19,19})) = \mathrm{Re}(-P)$ & $\mathrm{Im}(v_{19,19}) =     0$ & $\mathrm{Im}(i^n_{19,19}) = 0.243$ \\
        36.5 & $\mathrm{Re}(v_{19,20} \bar{i}^n_{19,20}) = \mathrm{Re}(P)$ & $\mathrm{Im}(v_{19,20}) = -0.243$ & $\mathrm{Im}(i^n_{19,20}) = -0.243$ \\
	\end{tabular}
	\end{small}
	\end{center}
\end{table}

\vspace{-1cm}
\begin{table}[H]
	\begin{center}
	\caption{Sorted $\alpha$ and $\Sigma$ values with corresponding equations and variables, respectively.}
	\begin{small}
	\begin{tabular}{c c c}
	\begin{tabular}{c||l}
		$\alpha_i$ & equation $i$ \\
	\hline
    4.69e+03 & $\abs{v_{5,1}} = V_g$ \\
        42.3 & $\abs{v_{20,20}} = V_g$ \\
        34.2 & $\abs{v_{18,20}} = V_g$ \\
          34 & $\abs{v_{19,19}} = V_g$ \\
        33.8 & $\abs{v_{20,18}} = V_g$ \\
        28.1 & $\abs{v_{16,20}} = V_g$ \\
	\end{tabular} & &
	\begin{tabular}{c||l}
		 $\sigma_{j,j}$ & variable $j$ \\
	\hline
   -1.04e+04 & $\mathrm{Im}(i^n_{5,1}) = 0.00128$ \\
   -1.02e+04 & $\mathrm{Im}(v_{5,1}) = 0.00197$ \\
    4.21e+03 & $\mathrm{Re}(v_{5,1}) = 0.0101$ \\
    2.12e+03 & $\mathrm{Re}(i^n_{5,1}) = -0.0101$ \\
        -279 & $\mathrm{Im}(v_{1,20}) = -0.243$ \\
        -212 & $\mathrm{Im}(v_{20,1}) = -0.243$ \\
	\end{tabular}
	\end{tabular}
	\end{small}
	\end{center}
\end{table}

\vspace{-.5cm}
WARNING: Newton-Raphson algorithm did not converge!

\newpage
Example 5: Changing some guess values (see table below).

\vspace{-0.5cm}
\begin{table}[H]
	\begin{center}
	\caption{Critical initial guess values compared to exact solution.}
	\begin{small}
	\begin{tabular}{l||c|c}
		variable & initial guess & exact solution \\
	\hline
		$v_{ 5, 1}$ & $  0.0101 + i\;   0.0020$ & $  1.0117 + i\;   0.1974$ \\
		$i^n_{ 5, 1}$  & $ -0.5067 + i\;   0.0638$ & $ -1.0133 + i\;   0.1276$ \\
	\end{tabular}
	\end{small}
	\end{center}
\end{table}

\vspace{-1cm}
\begin{table}[H]
	\begin{center}
	\caption{Sorted $\Gamma$ values with corresponding equation and variables.}
	\begin{small}
	\begin{tabular}{c||l|l|l}
		$\Gamma_{ijk}$ & equation $i$ & variable $j$ & variable $k$\\
	\hline
         101 & $\abs{v_{5,1}} = V_g$ & $\mathrm{Im}(v_{5,1}) = 0.00197$ & $\mathrm{Im}(v_{5,1}) = 0.00197$ \\
        10.4 & $\abs{v_{5,1}} = V_g$ & $\mathrm{Re}(v_{5,1}) = 0.0101$ & $\mathrm{Im}(v_{5,1}) = 0.00197$ \\
        2.46 & $\abs{v_{20,20}} = V_g$ & $\mathrm{Im}(v_{20,20}) =     0$ & $\mathrm{Im}(v_{20,20}) =     0$ \\
           2 & $\mathrm{Re}(v_{19,19} \bar{i}^n_{19,19})) = \mathrm{Re}(-P)$ & $\mathrm{Im}(v_{19,19}) =     0$ & $\mathrm{Im}(i^n_{19,19}) = 0.243$ \\
        1.83 & $\mathrm{Re}(v_{19,20} \bar{i}^n_{19,20}) = \mathrm{Re}(P)$ & $\mathrm{Im}(v_{19,20}) = -0.243$ & $\mathrm{Im}(i^n_{19,20}) = -0.243$ \\
        1.83 & $\mathrm{Re}(v_{20,19} \bar{i}^n_{20,19}) = \mathrm{Re}(P)$ & $\mathrm{Im}(v_{20,19}) = -0.243$ & $\mathrm{Im}(i^n_{20,19}) = -0.243$ \\
	\end{tabular}
	\end{small}
	\end{center}
\end{table}

\vspace{-1cm}
\begin{table}[H]
	\begin{center}
	\caption{Sorted $\alpha$ and $\Sigma$ values with corresponding equations and variables, respectively.}
	\begin{small}
	\begin{tabular}{c c c}
	\begin{tabular}{c||l}
		$\alpha_i$ & equation $i$ \\
	\hline
         123 & $\abs{v_{5,1}} = V_g$ \\
        1.61 & $\abs{v_{20,20}} = V_g$ \\
        1.08 & $\abs{v_{18,20}} = V_g$ \\
        1.06 & $\abs{v_{20,18}} = V_g$ \\
        1.04 & $\abs{v_{19,19}} = V_g$ \\
       0.774 & $\abs{v_{16,20}} = V_g$ \\
	\end{tabular} & &
	\begin{tabular}{c||l}
		 $\sigma_{j,j}$ & variable $j$ \\
	\hline
    1.07e+03 & $\mathrm{Im}(v_{5,1}) = 0.00197$ \\
       -20.5 & $\mathrm{Im}(i^n_{5,1}) = 0.0638$ \\
       -18.2 & $\mathrm{Re}(v_{5,1}) = 0.0101$ \\
        11.9 & $\mathrm{Re}(i^n_{5,1}) = -0.507$ \\
         5.8 & $\mathrm{Im}(v_{1,20}) = -0.243$ \\
        4.29 & $\mathrm{Im}(v_{20,1}) = -0.243$ \\
	\end{tabular}
	\end{tabular}
	\end{small}
	\end{center}
\end{table}

\vspace{-.5cm}
WARNING: Newton-Raphson algorithm did not converge!

\newpage
Example 6: Changing some guess values (see table below).

\vspace{-0.5cm}
\begin{table}[H]
	\begin{center}
	\caption{Critical initial guess values compared to exact solution.}
	\begin{small}
	\begin{tabular}{l||c|c}
		variable & initial guess & exact solution \\
	\hline
		$v_{ 5, 1}$ & $  0.0101 + i\;   0.0020$ & $  1.0117 + i\;   0.1974$ \\
	\end{tabular}
	\end{small}
	\end{center}
\end{table}

\vspace{-1cm}
\begin{table}[H]
	\begin{center}
	\caption{Sorted $\Gamma$ values with corresponding equation and variables.}
	\begin{small}
	\begin{tabular}{c||l|l|l}
		$\Gamma_{ijk}$ & equation $i$ & variable $j$ & variable $k$\\
	\hline
        1.43 & $\abs{v_{5,1}} = V_g$ & $\mathrm{Re}(v_{5,1}) = 0.0101$ & $\mathrm{Re}(v_{5,1}) = 0.0101$ \\
       0.998 & $\abs{v_{20,20}} = V_g$ & $\mathrm{Im}(v_{20,20}) =     0$ & $\mathrm{Im}(v_{20,20}) =     0$ \\
       0.602 & $\mathrm{Re}(v_{19,19} \bar{i}^n_{19,19})) = \mathrm{Re}(-P)$ & $\mathrm{Im}(v_{19,19}) =     0$ & $\mathrm{Im}(i^n_{19,19}) = 0.243$ \\
       0.578 & $\mathrm{Re}(v_{19,20} \bar{i}^n_{19,20}) = \mathrm{Re}(P)$ & $\mathrm{Im}(v_{19,20}) = -0.243$ & $\mathrm{Im}(i^n_{19,20}) = -0.243$ \\
       0.577 & $\mathrm{Re}(v_{20,19} \bar{i}^n_{20,19}) = \mathrm{Re}(P)$ & $\mathrm{Im}(v_{20,19}) = -0.243$ & $\mathrm{Im}(i^n_{20,19}) = -0.243$ \\
       0.574 & $\mathrm{Re}(v_{20,20} \bar{i}^n_{20,20})) = \mathrm{Re}(-P)$ & $\mathrm{Im}(v_{20,20}) =     0$ & $\mathrm{Im}(i^n_{20,20}) = 0.243$ \\
	\end{tabular}
	\end{small}
	\end{center}
\end{table}

\vspace{-1cm}
\begin{table}[H]
	\begin{center}
	\caption{Sorted $\alpha$ and $\Sigma$ values with corresponding equations and variables, respectively.}
	\begin{small}
	\begin{tabular}{c c c}
	\begin{tabular}{c||l}
		$\alpha_i$ & equation $i$ \\
	\hline
       0.996 & $\abs{v_{5,1}} = V_g$ \\
       0.307 & $\abs{v_{20,20}} = V_g$ \\
       0.128 & $\abs{v_{18,20}} = V_g$ \\
       0.127 & $\abs{v_{20,18}} = V_g$ \\
       0.112 & $\abs{v_{19,19}} = V_g$ \\
      0.0653 & $\abs{v_{18,18}} = V_g$ \\
	\end{tabular} & &
	\begin{tabular}{c||l}
		 $\sigma_{j,j}$ & variable $j$ \\
	\hline
          37 & $\mathrm{Im}(v_{5,1}) = 0.00197$ \\
        2.67 & $\mathrm{Re}(i^n_{5,1}) = -0.97$ \\
       -1.55 & $\mathrm{Re}(v_{5,1}) = 0.0101$ \\
       -1.54 & $\mathrm{Im}(v_{20,1}) = -0.243$ \\
        1.14 & $\mathrm{Im}(v_{20,19}) = -0.243$ \\
        1.14 & $\mathrm{Im}(v_{19,20}) = -0.243$ \\
	\end{tabular}
	\end{tabular}
	\end{small}
	\end{center}
\end{table}

\vspace{-.5cm}
SUCCESS: Newton-Raphson algorithm converged in 7 iterations!

\newpage
Example 7: Changing some guess values (see table below).

\vspace{-0.5cm}
\begin{table}[H]
	\begin{center}
	\caption{Critical initial guess values compared to exact solution.}
	\begin{small}
	\begin{tabular}{l||c|c}
		variable & initial guess & exact solution \\
	\hline
		$v_{ 5, 1}$ & $  0.5058 + i\;   0.0987$ & $  1.0117 + i\;   0.1974$ \\
		$v_{13,17}$ & $  0.1030 + i\;   0.0044$ & $  1.0298 + i\;   0.0439$ \\
		$v_{13,18}$ & $  0.0980 - i\;   0.0197$ & $  0.9802 - i\;   0.1973$ \\
		$v_{14,18}$ & $  0.1028 + i\;   0.0073$ & $  1.0282 + i\;   0.0731$ \\
		$v_{12,20}$ & $  0.1027 + i\;   0.0087$ & $  1.0271 + i\;   0.0870$ \\
		$v_{13,20}$ & $  0.0940 - i\;   0.0243$ & $  0.9402 - i\;   0.2435$ \\
		$i^n_{ 5, 1}$  & $ -0.1013 + i\;   0.0128$ & $ -1.0133 + i\;   0.1276$ \\
		$i^n_{13,17}$  & $ -0.0980 + i\;   0.0202$ & $ -0.9796 + i\;   0.2018$ \\
		$i^n_{13,18}$  & $  0.0980 - i\;   0.0197$ & $  0.9805 - i\;   0.1973$ \\
		$i^n_{14,18}$  & $ -0.0985 + i\;   0.0175$ & $ -0.9850 + i\;   0.1749$ \\
		$i^n_{13,20}$  & $  0.0997 - i\;   0.0258$ & $  0.9967 - i\;   0.2581$ \\
	\end{tabular}
	\end{small}
	\end{center}
\end{table}

\vspace{-1cm}
\begin{table}[H]
	\begin{center}
	\caption{Sorted $\Gamma$ values with corresponding equation and variables.}
	\begin{small}
	\begin{tabular}{c||l|l|l}
		$\Gamma_{ijk}$ & equation $i$ & variable $j$ & variable $k$\\
	\hline
        5.41 & $\mathrm{Re}(v_{20,1} \bar{i}^n_{20,1}) = \mathrm{Re}(P)$ & $\mathrm{Im}(v_{20,1}) = -0.243$ & $\mathrm{Im}(i^n_{20,1}) = -0.243$ \\
        4.83 & $\mathrm{Re}(v_{19,1} \bar{i}^n_{19,1})) = \mathrm{Re}(-P)$ & $\mathrm{Im}(v_{19,1}) =     0$ & $\mathrm{Im}(i^n_{19,1}) = 0.243$ \\
         4.8 & $\mathrm{Re}(v_{20,2} \bar{i}^n_{20,2})) = \mathrm{Re}(-P)$ & $\mathrm{Im}(v_{20,2}) =     0$ & $\mathrm{Im}(i^n_{20,2}) = 0.243$ \\
        4.15 & $\abs{v_{19,1}} = V_g$ & $\mathrm{Im}(v_{19,1}) =     0$ & $\mathrm{Im}(v_{19,1}) =     0$ \\
        4.12 & $\abs{v_{20,2}} = V_g$ & $\mathrm{Im}(v_{20,2}) =     0$ & $\mathrm{Im}(v_{20,2}) =     0$ \\
         3.8 & $\abs{v_{12,20}} = V_g$ & $\mathrm{Im}(v_{12,20}) = 0.0087$ & $\mathrm{Im}(v_{12,20}) = 0.0087$ \\
	\end{tabular}
	\end{small}
	\end{center}
\end{table}

\vspace{-1cm}
\begin{table}[H]
	\begin{center}
	\caption{Sorted $\alpha$ and $\Sigma$ values with corresponding equations and variables, respectively.}
	\begin{small}
	\begin{tabular}{c c c}
	\begin{tabular}{c||l}
		$\alpha_i$ & equation $i$ \\
	\hline
        3.49 & $\abs{v_{12,20}} = V_g$ \\
         3.4 & $\abs{v_{19,1}} = V_g$ \\
        3.37 & $\abs{v_{20,2}} = V_g$ \\
        2.48 & $\abs{v_{18,2}} = V_g$ \\
        2.43 & $\abs{v_{19,3}} = V_g$ \\
        2.38 & $\abs{v_{17,1}} = V_g$ \\
	\end{tabular} & &
	\begin{tabular}{c||l}
		 $\sigma_{j,j}$ & variable $j$ \\
	\hline
        -174 & $\mathrm{Im}(v_{12,20}) = 0.0087$ \\
        46.3 & $\mathrm{Re}(v_{13,20}) = 0.094$ \\
       -33.7 & $\mathrm{Im}(v_{19,1}) =     0$ \\
       -33.3 & $\mathrm{Im}(v_{20,2}) =     0$ \\
       -29.4 & $\mathrm{Im}(i^n_{19,1}) = 0.243$ \\
         -29 & $\mathrm{Im}(i^n_{20,2}) = 0.243$ \\
	\end{tabular}
	\end{tabular}
	\end{small}
	\end{center}
\end{table}

\vspace{-.5cm}
WARNING: Newton-Raphson algorithm did not converge!

\newpage
Example 8: Changing some guess values (see table below).

\vspace{-0.5cm}
\begin{table}[H]
	\begin{center}
	\caption{Critical initial guess values compared to exact solution.}
	\begin{small}
	\begin{tabular}{l||c|c}
		variable & initial guess & exact solution \\
	\hline
		$v_{ 5, 1}$ & $  0.5058 + i\;   0.0987$ & $  1.0117 + i\;   0.1974$ \\
		$v_{13,17}$ & $  0.1030 + i\;   0.0044$ & $  1.0298 + i\;   0.0439$ \\
		$v_{13,18}$ & $  0.0980 - i\;   0.0197$ & $  0.9802 - i\;   0.1973$ \\
		$v_{14,18}$ & $  0.1028 + i\;   0.0073$ & $  1.0282 + i\;   0.0731$ \\
		$v_{12,20}$ & $  0.8217 + i\;   0.0696$ & $  1.0271 + i\;   0.0870$ \\
		$v_{13,20}$ & $  0.0940 - i\;   0.0243$ & $  0.9402 - i\;   0.2435$ \\
		$i^n_{ 5, 1}$  & $ -0.1013 + i\;   0.0128$ & $ -1.0133 + i\;   0.1276$ \\
		$i^n_{13,17}$  & $ -0.0980 + i\;   0.0202$ & $ -0.9796 + i\;   0.2018$ \\
		$i^n_{13,18}$  & $  0.0980 - i\;   0.0197$ & $  0.9805 - i\;   0.1973$ \\
		$i^n_{14,18}$  & $ -0.0985 + i\;   0.0175$ & $ -0.9850 + i\;   0.1749$ \\
		$i^n_{13,20}$  & $  0.0997 - i\;   0.0258$ & $  0.9967 - i\;   0.2581$ \\
	\end{tabular}
	\end{small}
	\end{center}
\end{table}

\vspace{-1cm}
\begin{table}[H]
	\begin{center}
	\caption{Sorted $\Gamma$ values with corresponding equation and variables.}
	\begin{small}
	\begin{tabular}{c||l|l|l}
		$\Gamma_{ijk}$ & equation $i$ & variable $j$ & variable $k$\\
	\hline
        3.12 & $\mathrm{Re}(v_{13,20} \bar{i}^n_{13,20}) = \mathrm{Re}(P)$ & $\mathrm{Im}(v_{13,20}) = -0.0243$ & $\mathrm{Im}(i^n_{13,20}) = -0.0258$ \\
        2.85 & $\mathrm{Im}(v_{13,20} \bar{i}^n_{13,20}) = \mathrm{Im}(P)$ & $\mathrm{Im}(v_{13,20}) = -0.0243$ & $\mathrm{Re}(i^n_{13,20}) = 0.0997$ \\
         2.6 & $\abs{v_{13,17}} = V_g$ & $\mathrm{Im}(v_{13,17}) = 0.00439$ & $\mathrm{Im}(v_{13,17}) = 0.00439$ \\
        1.68 & $\abs{v_{14,18}} = V_g$ & $\mathrm{Im}(v_{14,18}) = 0.00731$ & $\mathrm{Im}(v_{14,18}) = 0.00731$ \\
        1.14 & $\mathrm{Re}(v_{12,20} \bar{i}^n_{12,20})) = \mathrm{Re}(-P)$ & $\mathrm{Im}(v_{12,20}) = 0.0696$ & $\mathrm{Im}(i^n_{12,20}) = 0.243$ \\
        1.14 & $\mathrm{Im}(v_{13,18} \bar{i}^n_{13,18}) = \mathrm{Im}(P)$ & $\mathrm{Im}(v_{13,18}) = -0.0197$ & $\mathrm{Re}(i^n_{13,18}) = 0.098$ \\
	\end{tabular}
	\end{small}
	\end{center}
\end{table}

\vspace{-1cm}
\begin{table}[H]
	\begin{center}
	\caption{Sorted $\alpha$ and $\Sigma$ values with corresponding equations and variables, respectively.}
	\begin{small}
	\begin{tabular}{c c c}
	\begin{tabular}{c||l}
		$\alpha_i$ & equation $i$ \\
	\hline
        2.37 & $\abs{v_{13,17}} = V_g$ \\
        1.47 & $\abs{v_{14,18}} = V_g$ \\
       0.748 & $\abs{v_{12,20}} = V_g$ \\
       0.265 & $\abs{v_{13,19}} = V_g$ \\
       0.242 & $\abs{v_{5,1}} = V_g$ \\
       0.222 & $\abs{v_{14,20}} = V_g$ \\
	\end{tabular} & &
	\begin{tabular}{c||l}
		 $\sigma_{j,j}$ & variable $j$ \\
	\hline
          72 & $\mathrm{Re}(v_{13,20}) = 0.094$ \\
         -66 & $\mathrm{Im}(v_{13,20}) = -0.0243$ \\
       -45.3 & $\mathrm{Im}(i^n_{13,20}) = -0.0258$ \\
        26.3 & $\mathrm{Re}(v_{13,18}) = 0.098$ \\
        26.3 & $\mathrm{Re}(i^n_{13,20}) = 0.0997$ \\
       -17.6 & $\mathrm{Im}(v_{13,18}) = -0.0197$ \\
	\end{tabular}
	\end{tabular}
	\end{small}
	\end{center}
\end{table}

\vspace{-.5cm}
WARNING: Newton-Raphson algorithm did not converge!

\newpage
Example 9: Changing some guess values (see table below).

\vspace{-0.5cm}
\begin{table}[H]
	\begin{center}
	\caption{Critical initial guess values compared to exact solution.}
	\begin{small}
	\begin{tabular}{l||c|c}
		variable & initial guess & exact solution \\
	\hline
		$v_{ 5, 1}$ & $  0.5058 + i\;   0.0987$ & $  1.0117 + i\;   0.1974$ \\
		$v_{13,17}$ & $  0.1030 + i\;   0.0044$ & $  1.0298 + i\;   0.0439$ \\
		$v_{13,18}$ & $  0.0980 - i\;   0.0197$ & $  0.9802 - i\;   0.1973$ \\
		$v_{14,18}$ & $  0.1028 + i\;   0.0073$ & $  1.0282 + i\;   0.0731$ \\
		$v_{12,20}$ & $  0.8217 + i\;   0.0696$ & $  1.0271 + i\;   0.0870$ \\
		$v_{13,20}$ & $  0.7522 - i\;   0.1948$ & $  0.9402 - i\;   0.2435$ \\
		$i^n_{ 5, 1}$  & $ -0.1013 + i\;   0.0128$ & $ -1.0133 + i\;   0.1276$ \\
		$i^n_{13,17}$  & $ -0.0980 + i\;   0.0202$ & $ -0.9796 + i\;   0.2018$ \\
		$i^n_{13,18}$  & $  0.0980 - i\;   0.0197$ & $  0.9805 - i\;   0.1973$ \\
		$i^n_{14,18}$  & $ -0.0985 + i\;   0.0175$ & $ -0.9850 + i\;   0.1749$ \\
		$i^n_{13,20}$  & $  0.0997 - i\;   0.0258$ & $  0.9967 - i\;   0.2581$ \\
	\end{tabular}
	\end{small}
	\end{center}
\end{table}

\vspace{-1cm}
\begin{table}[H]
	\begin{center}
	\caption{Sorted $\Gamma$ values with corresponding equation and variables.}
	\begin{small}
	\begin{tabular}{c||l|l|l}
		$\Gamma_{ijk}$ & equation $i$ & variable $j$ & variable $k$\\
	\hline
        4.62 & $\abs{v_{14,18}} = V_g$ & $\mathrm{Im}(v_{14,18}) = 0.00731$ & $\mathrm{Im}(v_{14,18}) = 0.00731$ \\
         4.4 & $\mathrm{Re}(v_{20,1} \bar{i}^n_{20,1}) = \mathrm{Re}(P)$ & $\mathrm{Im}(v_{20,1}) = -0.243$ & $\mathrm{Im}(i^n_{20,1}) = -0.243$ \\
         3.9 & $\mathrm{Re}(v_{19,1} \bar{i}^n_{19,1})) = \mathrm{Re}(-P)$ & $\mathrm{Im}(v_{19,1}) =     0$ & $\mathrm{Im}(i^n_{19,1}) = 0.243$ \\
        3.86 & $\mathrm{Re}(v_{20,2} \bar{i}^n_{20,2})) = \mathrm{Re}(-P)$ & $\mathrm{Im}(v_{20,2}) =     0$ & $\mathrm{Im}(i^n_{20,2}) = 0.243$ \\
         3.6 & $\abs{v_{13,17}} = V_g$ & $\mathrm{Im}(v_{13,17}) = 0.00439$ & $\mathrm{Im}(v_{13,17}) = 0.00439$ \\
        3.34 & $\abs{v_{19,1}} = V_g$ & $\mathrm{Im}(v_{19,1}) =     0$ & $\mathrm{Im}(v_{19,1}) =     0$ \\
	\end{tabular}
	\end{small}
	\end{center}
\end{table}

\vspace{-1cm}
\begin{table}[H]
	\begin{center}
	\caption{Sorted $\alpha$ and $\Sigma$ values with corresponding equations and variables, respectively.}
	\begin{small}
	\begin{tabular}{c c c}
	\begin{tabular}{c||l}
		$\alpha_i$ & equation $i$ \\
	\hline
        4.25 & $\abs{v_{14,18}} = V_g$ \\
        3.33 & $\abs{v_{13,17}} = V_g$ \\
        2.61 & $\abs{v_{19,1}} = V_g$ \\
        2.58 & $\abs{v_{20,2}} = V_g$ \\
        1.89 & $\abs{v_{18,2}} = V_g$ \\
        1.84 & $\abs{v_{19,3}} = V_g$ \\
	\end{tabular} & &
	\begin{tabular}{c||l}
		 $\sigma_{j,j}$ & variable $j$ \\
	\hline
       -34.4 & $\mathrm{Im}(v_{20,1}) = -0.243$ \\
        25.5 & $\mathrm{Im}(v_{19,1}) =     0$ \\
        25.3 & $\mathrm{Im}(v_{20,2}) =     0$ \\
       -25.1 & $\mathrm{Im}(i^n_{20,1}) = -0.243$ \\
       -23.8 & $\mathrm{Im}(v_{1,20}) = -0.243$ \\
        22.7 & $\mathrm{Re}(v_{13,18}) = 0.098$ \\
	\end{tabular}
	\end{tabular}
	\end{small}
	\end{center}
\end{table}

\vspace{-.5cm}
WARNING: Newton-Raphson algorithm did not converge!

\newpage
Example 10: Changing some guess values (see table below).

\vspace{-0.5cm}
\begin{table}[H]
	\begin{center}
	\caption{Critical initial guess values compared to exact solution.}
	\begin{small}
	\begin{tabular}{l||c|c}
		variable & initial guess & exact solution \\
	\hline
		$v_{ 5, 1}$ & $  0.5058 + i\;   0.0987$ & $  1.0117 + i\;   0.1974$ \\
		$v_{13,17}$ & $  0.1030 + i\;   0.0044$ & $  1.0298 + i\;   0.0439$ \\
		$v_{13,18}$ & $  0.0980 - i\;   0.0197$ & $  0.9802 - i\;   0.1973$ \\
		$v_{14,18}$ & $  0.8225 + i\;   0.0585$ & $  1.0282 + i\;   0.0731$ \\
		$v_{12,20}$ & $  0.8217 + i\;   0.0696$ & $  1.0271 + i\;   0.0870$ \\
		$v_{13,20}$ & $  0.7522 - i\;   0.1948$ & $  0.9402 - i\;   0.2435$ \\
		$i^n_{ 5, 1}$  & $ -0.1013 + i\;   0.0128$ & $ -1.0133 + i\;   0.1276$ \\
		$i^n_{13,17}$  & $ -0.0980 + i\;   0.0202$ & $ -0.9796 + i\;   0.2018$ \\
		$i^n_{13,18}$  & $  0.0980 - i\;   0.0197$ & $  0.9805 - i\;   0.1973$ \\
		$i^n_{14,18}$  & $ -0.0985 + i\;   0.0175$ & $ -0.9850 + i\;   0.1749$ \\
		$i^n_{13,20}$  & $  0.0997 - i\;   0.0258$ & $  0.9967 - i\;   0.2581$ \\
	\end{tabular}
	\end{small}
	\end{center}
\end{table}

\vspace{-1cm}
\begin{table}[H]
	\begin{center}
	\caption{Sorted $\Gamma$ values with corresponding equation and variables.}
	\begin{small}
	\begin{tabular}{c||l|l|l}
		$\Gamma_{ijk}$ & equation $i$ & variable $j$ & variable $k$\\
	\hline
        1.79 & $\abs{v_{13,17}} = V_g$ & $\mathrm{Im}(v_{13,17}) = 0.00439$ & $\mathrm{Im}(v_{13,17}) = 0.00439$ \\
        1.23 & $\mathrm{Im}(v_{13,18} \bar{i}^n_{13,18}) = \mathrm{Im}(P)$ & $\mathrm{Im}(v_{13,18}) = -0.0197$ & $\mathrm{Re}(i^n_{13,18}) = 0.098$ \\
       0.602 & $\mathrm{Re}(v_{13,18} \bar{i}^n_{13,18}) = \mathrm{Re}(P)$ & $\mathrm{Im}(v_{13,18}) = -0.0197$ & $\mathrm{Im}(i^n_{13,18}) = -0.0197$ \\
       0.425 & $\mathrm{Re}(v_{13,17} \bar{i}^n_{13,17})) = \mathrm{Re}(-P)$ & $\mathrm{Re}(v_{13,17}) = 0.103$ & $\mathrm{Re}(i^n_{13,17}) = -0.098$ \\
       0.207 & $\abs{v_{1,1}} = V_g$ & $\mathrm{Im}(v_{1,1}) =     0$ & $\mathrm{Im}(v_{1,1}) =     0$ \\
       0.195 & $\abs{v_{5,1}} = V_g$ & $\mathrm{Im}(v_{5,1}) = 0.0987$ & $\mathrm{Im}(v_{5,1}) = 0.0987$ \\
	\end{tabular}
	\end{small}
	\end{center}
\end{table}

\vspace{-1cm}
\begin{table}[H]
	\begin{center}
	\caption{Sorted $\alpha$ and $\Sigma$ values with corresponding equations and variables, respectively.}
	\begin{small}
	\begin{tabular}{c c c}
	\begin{tabular}{c||l}
		$\alpha_i$ & equation $i$ \\
	\hline
        1.61 & $\abs{v_{13,17}} = V_g$ \\
       0.112 & $\abs{v_{5,1}} = V_g$ \\
      0.0763 & $\abs{v_{1,1}} = V_g$ \\
      0.0381 & $\abs{v_{3,1}} = V_g$ \\
      0.0354 & $\abs{v_{12,20}} = V_g$ \\
      0.0342 & $\abs{v_{1,3}} = V_g$ \\
	\end{tabular} & &
	\begin{tabular}{c||l}
		 $\sigma_{j,j}$ & variable $j$ \\
	\hline
        26.5 & $\mathrm{Re}(v_{13,18}) = 0.098$ \\
         -18 & $\mathrm{Im}(v_{13,18}) = -0.0197$ \\
       -15.2 & $\mathrm{Im}(i^n_{13,18}) = -0.0197$ \\
       -9.01 & $\mathrm{Re}(i^n_{13,17}) = -0.098$ \\
        6.27 & $\mathrm{Im}(v_{13,17}) = 0.00439$ \\
       -3.71 & $\mathrm{Im}(v_{5,1}) = 0.0987$ \\
	\end{tabular}
	\end{tabular}
	\end{small}
	\end{center}
\end{table}

\vspace{-.5cm}
WARNING: Newton-Raphson algorithm did not converge!

\newpage
Example 11: Changing some guess values (see table below).

\vspace{-0.5cm}
\begin{table}[H]
	\begin{center}
	\caption{Critical initial guess values compared to exact solution.}
	\begin{small}
	\begin{tabular}{l||c|c}
		variable & initial guess & exact solution \\
	\hline
		$v_{ 5, 1}$ & $  0.5058 + i\;   0.0987$ & $  1.0117 + i\;   0.1974$ \\
		$v_{13,17}$ & $  0.8239 + i\;   0.0351$ & $  1.0298 + i\;   0.0439$ \\
		$v_{13,18}$ & $  0.0980 - i\;   0.0197$ & $  0.9802 - i\;   0.1973$ \\
		$v_{14,18}$ & $  0.8225 + i\;   0.0585$ & $  1.0282 + i\;   0.0731$ \\
		$v_{12,20}$ & $  0.8217 + i\;   0.0696$ & $  1.0271 + i\;   0.0870$ \\
		$v_{13,20}$ & $  0.7522 - i\;   0.1948$ & $  0.9402 - i\;   0.2435$ \\
		$i^n_{ 5, 1}$  & $ -0.1013 + i\;   0.0128$ & $ -1.0133 + i\;   0.1276$ \\
		$i^n_{13,17}$  & $ -0.0980 + i\;   0.0202$ & $ -0.9796 + i\;   0.2018$ \\
		$i^n_{13,18}$  & $  0.0980 - i\;   0.0197$ & $  0.9805 - i\;   0.1973$ \\
		$i^n_{14,18}$  & $ -0.0985 + i\;   0.0175$ & $ -0.9850 + i\;   0.1749$ \\
		$i^n_{13,20}$  & $  0.0997 - i\;   0.0258$ & $  0.9967 - i\;   0.2581$ \\
	\end{tabular}
	\end{small}
	\end{center}
\end{table}

\vspace{-1cm}
\begin{table}[H]
	\begin{center}
	\caption{Sorted $\Gamma$ values with corresponding equation and variables.}
	\begin{small}
	\begin{tabular}{c||l|l|l}
		$\Gamma_{ijk}$ & equation $i$ & variable $j$ & variable $k$\\
	\hline
        4.05 & $\mathrm{Re}(v_{20,1} \bar{i}^n_{20,1}) = \mathrm{Re}(P)$ & $\mathrm{Im}(v_{20,1}) = -0.243$ & $\mathrm{Im}(i^n_{20,1}) = -0.243$ \\
        3.96 & $\mathrm{Re}(v_{19,1} \bar{i}^n_{19,1})) = \mathrm{Re}(-P)$ & $\mathrm{Im}(v_{19,1}) =     0$ & $\mathrm{Im}(i^n_{19,1}) = 0.243$ \\
        3.92 & $\mathrm{Re}(v_{20,2} \bar{i}^n_{20,2})) = \mathrm{Re}(-P)$ & $\mathrm{Im}(v_{20,2}) =     0$ & $\mathrm{Im}(i^n_{20,2}) = 0.243$ \\
        3.54 & $\abs{v_{19,1}} = V_g$ & $\mathrm{Im}(v_{19,1}) =     0$ & $\mathrm{Im}(v_{19,1}) =     0$ \\
         3.5 & $\abs{v_{20,2}} = V_g$ & $\mathrm{Im}(v_{20,2}) =     0$ & $\mathrm{Im}(v_{20,2}) =     0$ \\
        2.96 & $\mathrm{Re}(v_{18,2} \bar{i}^n_{18,2})) = \mathrm{Re}(-P)$ & $\mathrm{Im}(v_{18,2}) =     0$ & $\mathrm{Im}(i^n_{18,2}) = 0.243$ \\
	\end{tabular}
	\end{small}
	\end{center}
\end{table}

\vspace{-1cm}
\begin{table}[H]
	\begin{center}
	\caption{Sorted $\alpha$ and $\Sigma$ values with corresponding equations and variables, respectively.}
	\begin{small}
	\begin{tabular}{c c c}
	\begin{tabular}{c||l}
		$\alpha_i$ & equation $i$ \\
	\hline
        2.76 & $\abs{v_{19,1}} = V_g$ \\
        2.72 & $\abs{v_{20,2}} = V_g$ \\
        2.17 & $\abs{v_{17,1}} = V_g$ \\
         2.1 & $\abs{v_{1,1}} = V_g$ \\
         2.1 & $\abs{v_{18,2}} = V_g$ \\
        2.03 & $\abs{v_{19,3}} = V_g$ \\
	\end{tabular} & &
	\begin{tabular}{c||l}
		 $\sigma_{j,j}$ & variable $j$ \\
	\hline
        29.5 & $\mathrm{Re}(v_{13,18}) = 0.098$ \\
        26.5 & $\mathrm{Im}(v_{20,1}) = -0.243$ \\
       -23.5 & $\mathrm{Im}(v_{13,18}) = -0.0197$ \\
       -20.6 & $\mathrm{Im}(i^n_{13,18}) = -0.0197$ \\
       -19.3 & $\mathrm{Im}(v_{19,1}) =     0$ \\
        19.2 & $\mathrm{Im}(i^n_{20,1}) = -0.243$ \\
	\end{tabular}
	\end{tabular}
	\end{small}
	\end{center}
\end{table}

\vspace{-.5cm}
SUCCESS: Newton-Raphson algorithm converged in 13 iterations!